\documentclass[onefignum,onetabnum]{siamart251216}

\usepackage{graphicx} 
\usepackage{algorithm}
\usepackage{algorithmic}

\usepackage[utf8]{inputenc}
\usepackage{graphicx}
\usepackage{bm}
\usepackage{comment}
\usepackage{latexsym,amsmath,amsfonts,amscd}
\usepackage{color}
\usepackage{bm}
\usepackage{tikz}
\usepackage{multirow}
\usepackage{booktabs}
\usepackage{xcolor}
\usepackage{hyperref}
\usepackage{setspace}
\usepackage{amssymb}
\usepackage{mathrsfs}

\usepackage[title]{appendix}
\usepackage{ulem}
\usepackage{tikz}

\usetikzlibrary{arrows,backgrounds,positioning}

\newtheorem{thm}{Theorem}[section]
\newtheorem{lem}{Lemma}[section]
\newtheorem{cor}[thm]{Corollary}

\newtheorem{rem}[thm]{Remark}
\newenvironment{myproof}{Proof:}{$\square$}

\newcommand{\bx}{{\bm{x}}}

\newcommand{\myitem}[1]{\noindent\underline{\bf #1}}
\newcommand{\beq}{\begin{equation}}
\newcommand{\eeq}{\end{equation}}

\newcommand{\bA}{{\bm{A}}}

\newcommand{\bpsi}{{\bm{\psi}}}

\newcommand{\mU}{{\mathcal{U}}}
\newcommand{\mV}{{\mathcal{V}}}
\newcommand{\mW}{{\mathcal{W}}}

\newcommand{\bmu}{\bm{\mu}}

\newcommand{\BOmega}{\boldsymbol{\Omega}}

\newcommand{\BF}{\mathbf{F}}

\usepackage{xcolor}

\newcommand{\pzc}[1]{{\color{black}#1}}
\newcommand{\ssy}[1]{{\color{black}#1}}

\title{Incremental Tensor-Train Compression from Streaming TT-Formatted Data with Applications to Reduced-Order Modeling of Linear Kinetic Equations}

\author{Shuyuan Shang\footnote{Department of Mathematics, The Hong Kong University of Science and Technology, Clear Water Bay, Kowloon, Hong Kong, China. E-mail: sshangad@connect.ust.hk.},\quad Wei Guo\footnote{Department of Mathematics and Statistics, Texas Tech University, Lubbock, TX, 70409, USA.},\quad   Zhichao Peng\footnote{Department of Mathematics, The Hong Kong University of Science and Technology, Clear Water Bay, Kowloon, Hong Kong, China. Corresponding author, E-mail: pengzhic@ust.hk.}}

\begin{document}

\maketitle

\begin{abstract}
\pzc{To mitigate the computational burden due to the high dimensionality of linear kinetic equations, low-rank solvers and data-driven reduced order models (ROMs) have been developed. Low-rank solvers based on matrix or tensor compression exploit low-rank structures within individual phase-space solutions, while data-driven ROMs exploit low-rank structures of the  solution manifold across parameters arising in  parametric multi-query applications such as sensitivity analysis and uncertainty quantification. It is attractive to exploit their complementary strengths by combining them. A key challenge is to extract a common reduced basis without reconstructing full solution snapshots or prescribing the compression ranks in advance. To address this challenge, we develop a deterministic incremental tensor-train (TT) compression algorithm that operates directly on streaming TT-formatted data and adapts the representation ranks to a prescribed approximation tolerance.} Given a new TT tensor, the proposed method updates an accumulated TT representation through core-wise projection, residual orthogonalization, and adaptive enrichment, retaining only the complementary information that cannot be represented within a prescribed tolerance. By operating entirely at the level of TT cores, the algorithm avoids reconstructing either the incoming tensor or the accumulated full tensor. We establish approximation error bounds for the proposed incremental approach. Moreover, we show that the accumulated TT representation provides a compressed analogue of standard proper orthogonal decomposition for full-order
snapshot data. Using this connection,  we construct reduced order models for kinetic equations directly from
streaming low-rank solution data through operations on TT cores, without first reconstructing full snapshots. Numerical experiments on parametric radiative transfer equations demonstrate
that the proposed method achieves comparable reconstruction accuracy with
substantially reduced wall time and yields efficient and accurate ROMs directly
from compressed low-rank data.
\end{abstract}


\section{Introduction\label{sec:introduction}}


\pzc{The high dimensionality of kinetic equations has motivated the development of dimensionality reduction techniques, including low-rank solvers
\cite{bachmayr2023low,einkemmer2025review} and data-driven reduced-order models (ROMs)
\cite{quarteroni2011certified,benner2015survey}.  In this paper, we focus on linear kinetic equations, also known as radiative transfer equation or Boltzmann transport equation, modeling photon and neutron transport. Low-rank solvers e.g. \cite{einkemmer2021asymptotic,ding2021dynamical,ceruti2022rank,peng2023sweep,einkemmer2024asymptotic,kusch2023robust,truong2024tensor,sands2025high,bachmayr2026low,guo2026inexact,guo2026highly,deshpande2026multigroup} exploit phase-space low-rank structure within individual solutions and compute these solutions directly in compressed matrix or tensor formats. Data-driven ROMs e.g. \cite{mcclarren2019calculating,buchan2015pod,choi2021space,tano2021affine,peng2022reduced,mcclarren2022data,behne2023parametric,coale2023reduced,halvic2023non,coale2024reduced,hardy2024proper,peng2024micro,matsuda2025reduced,pan2026data} target at parametric equations arising from multi-query applications such as sensitivity analysis and uncertainty quantification, and  exploit low-rank structures of the parametric solution manifold across parameters. For more complete literature review, we refer readers to the review papers for low-rank solvers \cite{bachmayr2023low,einkemmer2025review} and recent data-driven ROM \cite{pan2026data}. These complementary strengths motivate the construction of data-driven ROMs directly from low-rank solution data. A central challenge is to construct a representative reduced basis incrementally from such data without reconstructing full solution snapshots or prescribing the compression ranks in advance.}

Tensor decompositions provide a natural framework for representing and compressing solution data. Among them, the tensor-train (TT) format is particularly attractive because it
provides stable, quasi-optimal truncation while achieving storage and rounding
costs that scale linearly with the tensor dimension when
the TT ranks remain moderate
\cite{oseledets2011tensor,grasedyck2013literature}. In many applications, however, solution snapshots are generated sequentially, motivating streaming TT algorithms that update compressed representations as new samples become available. This has led to growing interest in streaming TT algorithms, which update compressed tensor representations as new samples become available. Streaming TT approximation (STTA) \cite{kressner2023streaming} constructs TT representations from randomized sketches, whereas deterministic approaches such as TT-FOA \cite{abed2021adaptive} and TT-ICE \cite{aksoy2024incremental} incrementally update a TT representation as new tensor data become available.

For solution data generated by low-rank solvers, a key limitation concerns the representation of the incoming data. The deterministic incremental methods discussed above assume access to incoming data through dense tensor entries, whereas low-rank solvers produce solutions directly in compressed matrix or tensor formats. Reconstructing full solution tensors solely to apply a streaming compression algorithm can incur storage and computational costs, undermining the computational advantages of the low-rank solver. This motivates a deterministic, rank-adaptive incremental compression method that operates directly on streaming TT-formatted data.

In addition, a similar challenges also arises when constructing ROMs directly from these compressed low-rank data. Classical ROMs \cite{benner2015survey} are constructed from full-order solution data. Recent tensor-based ROMs \cite{ballani2016reduced,mamonov2022interpolatory,mamonov2024tensorial,
mueller2026tensor,saibaba2025tensor,he2026tensor} further introduce tensor compression or tensorial operations
into the ROM construction, but their construction still relies on full-order data. Although tensor-based dynamic mode decomposition \cite{klus2018tensor} can
handle time-history data in TT format, it does not address the incremental arrival of individual snapshots already in low-rank form. 

In this work, to enable efficient compression and ROM construction from streaming low-rank data, we develop a deterministic incremental TT compression algorithm that operates directly on streaming TT-formatted data. Given a new TT tensor, the proposed method updates an accumulated TT representation through core-wise projection, residual orthogonalization, and adaptive enrichment. At each mode, the incoming TT core is projected onto the current TT subspace, and only information that cannot be represented within a prescribed tolerance is retained. As a result, the algorithm operates entirely on TT cores and avoids reconstructing either the incoming tensor or the full accumulated tensor. Beyond data compression, we show that the accumulated TT tensor can be
interpreted as a decomposition into a compressed reduced basis and
training-sample coefficients, enabling efficient model order reduction directly from low-rank data.

The main contributions of this work are summarized as follows:
\begin{itemize}
\item \textbf{Incremental TT construction from streaming TT-formatted data.}
We propose a deterministic algorithm that incrementally constructs an accumulated TT
representation directly from TT-formatted inputs without reconstructing either the incoming full tensors or the accumulated full
tensor. The method performs projection, residual orthogonalization, and
adaptive enrichment through operations on TT cores, selecting representation
ranks according to a prescribed approximation tolerance. We establish
reconstruction error bounds and derive complexity estimates showing linear
dependence on the number of tensor modes when mode sizes and TT ranks remain
bounded.

\item \textbf{Application to reduced-order modeling.}
We show that the accumulated TT representation provides a compression of
solution data analogous to proper orthogonal decomposition (POD), thereby
enabling efficient low-rank ROMs for kinetic equations directly construct from the output of modern
low-rank solvers. Numerical experiments demonstrate lower offline computational and storage costs than classical ROMs built from full-order
snapshots, and lower online computational costs than the underlying
low-rank solver.
\end{itemize}


The remainder of
the paper is organized as follows. Section~2 reviews background material on low-rank solvers for kinetic equations and
tensor-train representations. Section~3 presents the incremental TT compression
algorithm. Section~4 establishes the approximation properties and computational
complexity. Section~5 develops reduced-order modeling formulations based on the
accumulated TT representation. Section~6 reports numerical experiments, and
Section~7 concludes the paper.
\section{Low-rank solvers for kinetic equations and tensor train representations}
Kinetic equations provide mesoscopic models of interacting particle systems. In this paper, we focus on linear kinetic equations for photon or neutron transport. In multi-query applications such as sensitivity analysis, desing optimization, imaging and uncertainty quantification, these linear kinetic equations are often solved for parameter values determining physical properties such as material or source configurations. As a model problem for such parametric setting, we consider the following
steady-state linear kinetic equation:
\begin{subequations}
\label{eq:rte}
\begin{align}
    &\BOmega\cdot\nabla_\bx f(\bx,\BOmega;\bmu)
    +\sigma_t(\bx;\bmu)f(\bx,\BOmega;\bmu)
    =
    \sigma_s(\bx;\bmu)\rho(\bx;\bmu)+G(\bx;\bmu),\\
   &\rho(\bx;\bmu)
    = \frac{1}{4\pi}\int_{\mathbb{S}^2}f(\bx,\BOmega;\bmu)\,d\BOmega,
    \qquad
    \sigma_t(\bx;\bmu)=\sigma_s(\bx;\bmu)+\sigma_a(\bx;\bmu),
\end{align}
\end{subequations}
where \(f(\bx,\BOmega;\bmu)\) denotes the particle distribution at the spatial location  and the particle transport direction \(\BOmega\in\mathbb{S}^2\), with the parameter \(\bmu\) determining material or source configurations. Furthermore, \(\rho(\bx;\bmu)\) denotes the macroscopic density, and \(\sigma_s(\bx;\bmu)\), \(\sigma_a(\bx;\bmu)\), and
\(\sigma_t(\bx;\bmu)\) are the parametric scattering, absorption, and total cross sections, respectively. The parametric source term is denoted by \(G(\bx;\bmu)\). Homogeneous inflow boundary conditions are imposed  for simplicity.

The dependence of $f$ on both spatial and angular variables
makes fully resolved discretizations expensive, particularly
when solutions are required at many parameter values.  Recently, to mitigate this challenge, low-rank solvers  have been developed for dimensionality reduction through matrix or tensor compressions to exploit the phase-space low-rank structures of the solution. For example,  \cite{einkemmer2020asymptotic,haut2026efficient,guo2026highly} use a low-rank matrix
approximation in the form $\BF(\bmu)\approx X(\bmu)S(\bmu)V(\bmu)^T$, where
$F(\bmu)\in\mathbb{R}^{N_{\bx}\times N_{\BOmega}}$ denotes the matrix collecting discrete spatial degrees of freedom (DOFs) for each angle, 
$X(\bmu)\in\mathbb{R}^{N_{\bx}\times r}$, $V(\bmu)\in\mathbb{R}^{N_{\BOmega}\times r}$ and
$S(\bmu)\in\mathbb{R}^{r\times r}$
 are low-rank spatial and angular basis and coefficient core. Solvers in
 \cite{truong2024tensor,deshpande2026multigroup} use tensor-train (TT)
decompositions.

Data-driven reduced-order models (ROMs) have been actively developed for parametric linear kinetic equations by extracting low-rank structures across parameters from solution data for training parameters. 
However, they are typically built from full
solution snapshots, without directly exploiting the compressed
representations produced by low-rank solvers. To exploit these
benefits in ROM construction, we will adapt the tensor train framework to incrementally build 
reduced basis from streaming low-rank solution data. Before moving to our approach, we first briefly review the TT decomposition.
\subsection{Tensor train decomposition}
The TT decomposition is an efficient low-rank representation for high-dimensional tensors, effectively mitigating the curse of dimensionality. Given a \(d\)-dimensional tensor \(\bm{A} \in \mathbb{R}^{n_1 \times n_2 \times \cdots \times n_d}\), the TT format represents \(\bm{A}\) as a sequence of third-order tensors (TT cores) through a matrix product structure
\begin{equation}
\label{eq:tt}
    \bm{A}(i_1,\dots,i_d) = \mathcal{V}_1(1,i_1,:) \, \mathcal{V}_2(:,i_2,:) \cdots \mathcal{V}_d(:,i_d,1),
\end{equation}
where each core \(\mathcal{V}_i \in \mathbb{R}^{r_{i-1} \times n_i \times r_i}\), and \(\{r_i\}_{i=0}^d\) are the TT ranks with \(r_0 = r_d = 1\). The TT ranks quantify the coupling between adjacent modes and determine both the approximation accuracy and the computational cost. For brevity, we may suppress the indices in \eqref{eq:tt} and write
\(
\bm{A} = \mathcal{V}_1\,\mathcal{V}_2 \cdots \mathcal{V}_d.
\)

The storage complexity of the TT representation is
$
O\!\left(\sum_{i=1}^d n_i r_{i-1} r_i\right),
$
which, under the assumption \(n_i \sim n\) and \(r_i \sim r\), scales as \(O(d n r^2)\). This linear scaling in the dimension \(d\) contrasts sharply with the exponential complexity \(O(n^d)\) required to store the full tensor.

Given a full tensor \(\bm{B}\in\mathbb{R}^{n_1 \times n_2 \times \cdots \times n_d}\), the standard construction of its TT decomposition \(\bm{A}\) is the \texttt{TT-SVD} algorithm \cite{oseledets2011tensor}, which relies on successive tensor unfoldings (also called matricizations) and truncated SVDs. The mode-\(i\) unfolding reshapes \(\bm{B}\) into a matrix, denoted by
$$
B_{\langle i \rangle}=\texttt{reshape}(\bm{B},[n_1\cdots n_i, \;n_{i+1}\cdots n_{d}]) \in \mathbb{R}^{(n_1\cdots n_i) \times (n_{i+1} \cdots n_d)}.
$$
Moreover, a tensor is uniquely determined by any of its unfoldings; that is, two tensors are identical if and only if their mode-\(i\) unfoldings coincide. Applying a truncated SVD to the first unfolding,
\(
B_{\langle 1 \rangle} \approx U_1 S_1 V_1^\top,
\)
yields a rank-\(r_1\) approximation with error bounded by \(\epsilon_1\) in the Frobenius norm. The first core \(\mathcal{V}_1\) of \(\bm{A}\) is obtained by reshaping
\(U_1\in\mathbb{R}^{n_1\times r_1}\) into a tensor of size
\(1\times n_1\times r_1\). The remaining factor
\(S_1V_1^\top\in\mathbb{R}^{r_1\times(n_2\cdots n_d)}\) is interpreted
as a tensor of size \(r_1\times n_2\times\cdots\times n_d\). Grouping
the first two indices gives a matrix
\(
C_2\in\mathbb{R}^{(r_1n_2)\times(n_3\cdots n_d)},
\)
to which the next truncated SVD with rank $r_2$ and approximation error $\epsilon_2$,
\(
C_2\approx U_2S_2V_2^\top,
\)
is applied. The second core \(\mathcal{V}_2\) is then obtained by
reshaping \(U_2\in\mathbb{R}^{r_1n_2\times r_2}\) into a tensor of size \(r_1\times n_2\times r_2\). Repeating this procedure, at the \(i\)-th step, a truncated SVD produces factors \(U_i S_i V_i^\top\) with approximation error \(\epsilon_i\), from which the core \(\mathcal{V}_i \in \mathbb{R}^{r_{i-1} \times n_i \times r_i}\) is obtained by reshaping \(U_i\in\mathbb{R}^{r_{i-1}n_i\times r_i}\).This process continues until all cores are constructed. The computational complexity of \texttt{TT-SVD} is dominated by the SVDs of the unfolding matrices and scales as
{
\(
O\!\left(\sum_{i=1}^{d-1} r_{i-1} r_i n^{\,d-i+1}\right),
\)
}
assuming \(n_i \sim n\). The resulting approximation is quasi-optimal in the sense that the error is within a small constant factor of the best achievable TT approximation with the same ranks, with error bounded by
\(
\|\bm{B} - \bm{A}\|_F \;\le\; \sqrt{\sum_{i=1}^{d-1} \epsilon_i^2}.
\)

If a tensor is already given in TT format, it can be recompressed using the \texttt{TT-rounding} algorithm with a prescribed threshold $\epsilon_{\rm tt}$. This procedure consists of successive orthogonalizations (via QR factorizations) and truncated SVDs applied to matricizations of the TT cores, producing a compressed TT representation with reduced ranks. The resulting approximation is quasi-optimal, as \texttt{TT-rounding} is mathematically equivalent to \texttt{TT-SVD} applied to the full tensor. The computational complexity of \texttt{TT-rounding} scales as $O(d n r^3)$, thereby avoiding the exponential scaling with respect to \(d\).

\section{Incremental construction of TT reduced basis from low-rank data\label{sec:LR-iTT-F}}
To exploit the complementary strengths of low-rank solvers
and data-driven ROMs in capturing phase-space and
cross-parameter low-rank structures, respectively, we develop
an incremental algorithm, denoted by LR-iTT, for constructing
a TT reduced basis directly from streaming low-rank data.
The algorithm builds on the project-then-enrich strategy
used in incremental SVD \cite{brand2006fast} and incremental
TT compression of streaming full tensor data \cite{aksoy2024incremental}
to construct an accumulated TT representation.
To operate directly on compressed TT data, we perform
core-wise projection and residual enrichment without
reconstructing either the incoming snapshots or the
accumulated full tensor. We first present the construction
with full residual enrichment for informative incoming data, denoted by LR-iTT-F.
Then, to further accelerate the algorithm, we introduce
a reduced-enrichment variant, denoted by LR-iTT-R, which
selectively retains residual directions to limit rank growth
and reduce computational costs. Finally, we discuss the
relation of the proposed methods to existing incremental
TT algorithms.

\subsection{Incremental TT construction for low-rank data with full residual enrichment (LR-iTT-F)}
\label{sec:LR-iTT-F-intro}
Consider a sequence of low-rank data \(\{\bm{A}^k\}_{k=1}^{n_\mu}\), where each tensor \(\bm{A}^k \in \mathbb{R}^{n_1\times\cdots\times n_d}\) is already available in the TT format \eqref{eq:tt}. Specifically, each \(\bm{A}^k\) is characterized by TT cores \(\mathcal{V}_i^k \in \mathbb{R}^{s_{i-1}\times n_i \times s_i}\) with TT ranks \(\{s_i\}_{i=0}^d\) satisfying \(s_0 = s_d = 1\).
The proposed incremental algorithm progressively constructs a $(d+1)$-dimensional accumulated tensor $\bpsi\in \mathbb{R}^{n_1\times\cdots\times n_d\times n_\mu}$ in the TT format
\begin{equation}
\bpsi= \mU_1\,\mU_2\,\cdots \mU_{d}\,\mU_{d+1}
\end{equation}
with TT cores $\{\mU_i\}_{i=1}^{d+1}$ and TT ranks $\{r_i\}_{i=0}^{d+1}$, $r_0=r_{d+1}=1$,
such that \ssy{$\| \bm{A}^k-\bpsi(:,:,\cdots,:,k) \|_F\lesssim \sqrt{d} \epsilon$, $k=1,\ldots,n_\mu$, where $\epsilon$ is the prescribed per-core residual tolerance
applied at each of the physical modes.}

To initialize the algorithm, we may construct $\bpsi^1$ directly from the first tensor $\bA^1$ and start the incremental updates at $k=2$. Alternatively, warm-start initialization based on prior knowledge of the dataset is also possible. Suppose that after processing the first $k-1$ tensors, we have already constructed $\bpsi^{k-1}$ with TT cores $\{\mU_i^{k-1}\}_{i=1}^{d+1}$ as the current TT approximation of the accumulated tensor. When a new tensor $\bm{A}^{k}$ with cores $\{\mV^k_i\}_{i=1}^{d}$ becomes available, the TT cores of $\bpsi^{k-1}$ are updated to produce $\bpsi^{k}$ while maintaining the prescribed accuracy guarantee. Furthermore, we assume that each core of the accumulated TT \(\mathcal{U}_i^{k-1}\) is \textit{left-orthogonal}, meaning that its left unfolding
\[
U_i^{k-1} := \texttt{reshape}\!\left(\mathcal{U}_i^{k-1},\,[r_{i-1}^{k-1} n_i,\; r_i^{k-1}]\right)
\]
has orthonormal columns. We also assume that the cores $\{\mV^k_i\}_{i=1}^{d}$ of the incoming TT tensor $\bm{A}^{k}$ are \textit{right-orthogonal}. Specifically, for each $i$, the right unfolding
\[
V^k_i := \texttt{reshape}\!\left(\mathcal{V}^k_i,\,[s_{i-1},\; n_i s_i]\right),
\]
has orthonormal rows.

\myitem{(a) Processing the first mode.}

Define \(\mathcal{W}_1 := \mathcal{V}_1^{k}\) and 
\[
W_1 := \texttt{reshape}\!\left(\mathcal{W}_1,\,[n_1,\, s_1]\right).
\]
We project \(W_1\) onto the orthogonal complement of the column space of \(U_1^{k-1}\) and define the residual in mode-1
\[
R_1 = \big(I - U_1^{k-1} (U_1^{k-1})^T\big) W_1,
\]
which represents the component of the new tensor not captured by the current TT subspace. Applying a QR decomposition yields
\[
P_1 = \texttt{orth}(R_1),
\]
whose columns form an orthonormal basis for \(\mathrm{range}(R_1)\). By construction,
\begin{equation}
\big(U_1^{k-1} (U_1^{k-1})^T + P_1P_1^T\big) W_1 = W_1.
\label{eq:projection_1}
\end{equation}

\myitem{(b) Updating $\mU^{k-1}_i\rightarrow\mU^{k}_{i}$ and processing the next mode.} 

For \(i=1\), no padding is required, and we set
\(\widetilde{\mU}^{k-1}_1=\mU^{k-1}_1\) and
\(\widetilde{U}^{k-1}_1=U^{k-1}_1\).  For \(i>1\), as noted in
\cite{liu2018incremental,aksoy2024incremental}, the rank change in the previous
TT core update for \(\mU_{i-1}^{k}\) requires a zero-padding procedure for the
current core \(\mU^{k-1}_i\) to restore \textit{dimension consistency} of the
TT representation.  This procedure yields the padded core
\(\widetilde{\mU}^{k-1}_i\) of size
\(r_{i-1}^{k}\times n_i\times r_i^{k-1}\) and the associated left unfolding
\(\widetilde{U}^{k-1}_i\).  The padded core preserves the left-orthogonality of
\(\mU^{k-1}_i\).

If the residual satisfies \(\|R_i\|_F > \epsilon\), the mode-$i$ component of $\bm A^k$ contains significant new information and should be incorporated into the accumulated subspace $\mathrm{range}(\widetilde{U}_i^{k-1})$. In this case, the increment \(P_i\) is appended to \(\widetilde{U}_i^{k-1}\), yielding
\begin{equation}
\label{eq:Uupdate}
U_i^{k} = \big[\widetilde{U}_i^{k-1} \;\; P_i\big],
\end{equation}
whose columns remain orthonormal. This construction satisfies an exactness property (see Lemma~\ref{lem:error}). The rank is updated as \(r_i^{k} = r_i^{k-1} + p_i\), where $p_i \le s_i$ is the number of columns of $P_i$.
Otherwise, the new information is already well represented by $\widetilde{U}_i^{k-1}$ up to the prescribed tolerance, and no enrichment is performed. We set \(U_i^{k} = \widetilde{U}_i^{k-1}\) with \(r_i^{k} = r_i^{k-1}\), yielding a controlled approximation error bound (see Lemma~\ref{lem:error}). The updated core is obtained by reshaping
\[
\mathcal{U}_i^{k} = \texttt{reshape}\!\left(U_i^{k},\,[r_{i-1}^{k},\, n_i,\, r_i^{k}]\right).
\]

If \(i < d\), we proceed to the next mode and compute
\begin{equation}
\label{eq:W}
\mathcal{W}_{i+1} := \texttt{reshape}\!\left((U_i^{k})^T W_i V^k_{i+1},\,[r_i^{k},\, n_{i+1},\, s_{i+1}]\right),
\end{equation}
and its left unfolding
\(
W_{i+1} := \texttt{reshape}\!\left(\mathcal{W}_{i+1},\,[r_i^{k} n_{i+1},\, s_{i+1}]\right).
\)
The residual with respect to the padded matrix \(\widetilde{U}_{i+1}^{k-1}\) is defined as:
\[
R_{i+1} = \big(I - \widetilde{U}_{i+1}^{k-1} (\widetilde{U}_{i+1}^{k-1})^T\big) W_{i+1}.
\]
Applying a QR decomposition yields an orthonormal increment
\[
P_{i+1} = \texttt{orth}(R_{i+1}),
\]
ensuring the decomposition holds:
\begin{equation}
\big(\widetilde{U}_{i+1}^{k-1} (\widetilde{U}_{i+1}^{k-1})^T + P_{i+1} P_{i+1}^T\big) W_{i+1} = W_{i+1}.
\label{eq:projection_i_p1}
\end{equation}
The mode-one projection identity \eqref{eq:projection_1} is the special case of
\eqref{eq:projection_i_p1} obtained from
\(\widetilde{U}^{k-1}_{1}=U^{k-1}_1\).

\myitem{(c) Updating the last core.} 

The update of the last core is slightly different. The matrix \(W_d\) is projected onto the column space of \(U_d^{k}\), yielding the slice
\[
W_{d+1} = (U_d^{k})^T W_d \in \mathbb{R}^{r_d^{k} \times 1}.
\]
Applying the zero-padding procedure to the last core \(\mathcal{U}_{d+1}^{k-1}\) produces the padded core \(\widetilde{\mathcal{U}}_{d+1}^{k-1}\in\mathbb{R}^{r_d^k\times (k-1)\times 1}\). We then append $W_{d+1}$ as the $k$-th slice to form the updated core
$\mathcal{U}_{d+1}^k \in \mathbb{R}^{r_d^k \times k \times 1}$:
\[\mathcal{U}_{d+1}^{k}(:, 1:k-1, 1) = \widetilde{\mathcal{U}}_{d+1}^{k-1}, \quad \mathcal{U}_{d+1}^{k}(:, k, 1) = W_{d+1}.\]
Consequently, the updated accumulated tensor is given by
\[
\bpsi^{k}
=
\mU^{k}_1\mU^{k}_2\cdots \mU^{k}_{d}\,\mU^{k}_{d+1}.
\]

\myitem{(d) Re-orthogonalization and truncation for rank control.}

If the TT ranks of $\bpsi^{k}$ exceed a prescribed maximum $r_{\mathrm{max}}$, \texttt{TT-rounding} is performed to reorthogonalize the TT cores, remove redundant information, and reduce the TT ranks of \(\bpsi^{k}\) while maintaining the prescribed accuracy.

The proposed construction employs an incremental project-and-enrich strategy to build the reduced basis and, unlike existing approaches \cite{mueller2026tensor}, avoids the exponential scaling with the dimension. The overall procedure is summarized in Algorithm~\ref{alg:proposed_tt_update}.

\subsection{Further acceleration through reduced enrichment}

In the full enrichment approach discussed in Section~\ref{sec:LR-iTT-F-intro}, when the residual $R_i$ at
the $i$-th core satisfies $\|R_i\|_F>\epsilon$, we expand
the core basis by appending an orthonormal basis $P_i$ for
the entire residual space. However, this strategy can be
too aggressive. The new information contained in $R_i$ may
be well represented by only a few dominant directions.
Retaining the remaining directions may provide little
improvement in accuracy while increasing computational
and storage costs. 

To address this issue, we introduce a reduced enrichment
strategy that only adds selected dominant directions of the residual space, ensuring the orthonormal basis of their span $P_i$ satisfying $\|(I-P_iP_i^T)R_i\|_F\leq \epsilon$. The detailed selection procedure is as follows.
For the residual $R_i\in\mathbb{R}^{m\times s}$, define
\begin{equation}
    a_j=\|R_i(:,j)\|_2,\qquad j=1,\ldots,s,
    \qquad
    a_{\max}=\max_{1\leq j\leq s}a_j.
\end{equation}
Given a parameter $\alpha\in(0,1]$, we identify the weak
columns eligible for removal by
\begin{equation}
    L=\{j:a_j<\alpha a_{\max}\}.
\end{equation}
All columns outside $L$ are retained. We order the indices
in $L$ as $j_1,\ldots,j_{|L|}$ such that
$a_{j_1}\leq\cdots\leq a_{j_{|L|}}$, and determine the
number of unimportant directions to discard by
\begin{equation}
    q=\max\big\{
    \ell\in\{0,\ldots,|L|\}:
    \sum_{t=1}^{\ell}a_{j_t}^2\leq\epsilon^2
    \big\}.
\end{equation}
The empty sum is taken as zero, and $q=0$ means that
no columns are discarded. The enrichment matrix is then
constructed as
\begin{equation}
    P_i=\texttt{orth}(R_i(:,S)),
    \qquad
    S=\{1,\ldots,s\}\setminus\{j_1,\ldots,j_q\}.
\end{equation}
Since the retained columns lie in
$\operatorname{range}(P_i)$ and orthogonal projection
does not increase the Euclidean norm, we have
$\|(I-P_iP_i^T)R_i\|_F^2
    \leq\sum_{t=1}^{q}a_{j_t}^2
    \leq\epsilon^2.$

The parameter $\alpha$ controls which columns are protected
during selection. A smaller $\alpha$ requires retaining
more columns and makes the selection more conservative,
whereas a larger $\alpha$ makes more columns eligible
for removal. In all cases, the total squared norm of
the discarded columns is bounded by $\epsilon^2$.

We summarize the incremental TT compression algorithm
with full and reduced enrichment, denoted by LR-iTT-F
and LR-iTT-R, respectively, in
Algorithm~\ref{alg:proposed_tt_update}. The residual
basis selection procedure used by LR-iTT-R is detailed
in Algorithm~\ref{alg:reduced_basis}.

\begin{algorithm}[t]
\footnotesize
\caption{\ssy{Incremental update from streaming TT data}
\label{alg:proposed_tt_update}}
\begin{algorithmic}[1]
\REQUIRE Left-orthogonal TT cores \(\{\mU_i^{k-1}\}_{i=1}^{d+1}\) of
\(\bpsi^{k-1}\); right-orthogonal TT cores
\(\{\mV^k_i\}_{i=1}^{d}\) of the new snapshot \(\bm A^k\);
residual tolerance \(\epsilon\); TT rounding tolerance \(\epsilon_{\rm tt}\);
maximum rank \(r_{\max}\); enrichment rule \(\in\{\texttt{LR-iTT-F},\texttt{LR-iTT-R}\}\);
parameter \(\alpha\in(0,1]\) (LR-iTT-R only).
\ENSURE Updated TT cores \(\{\mU_i^{k}\}_{i=1}^{d+1}\) of \(\bpsi^k\).
\STATE \(\mW_1 \leftarrow \mV_1\),
\(W_1 \leftarrow \texttt{reshape}(\mW_1,[n_1,s_1])\) \hfill \(\triangleright\) {\(\mW_1\in\mathbb{R}^{1\times n_1\times s_1}\)}
\FOR{\(i=1\) to \(d\)}
    \IF{\(i=1\)}
        \STATE \(\widetilde{\mU}_1^{k-1}\leftarrow \mU_1^{k-1}\)
    \ELSE
        \STATE Zero-pad \(\mU_i^{k-1}\) to obtain
        \(\widetilde{\mU}_i^{k-1}\in
        \mathbb{R}^{r_{i-1}^k\times n_i\times r_i^{k-1}}\)
    \ENDIF
    \STATE \(\widetilde U_i^{k-1}\leftarrow
        \texttt{reshape}(\widetilde{\mU}_i^{k-1},
        [r_{i-1}^k n_i,r_i^{k-1}])\)
    \STATE \(R_i \leftarrow
    \bigl(I-\widetilde U_i^{k-1}(\widetilde U_i^{k-1})^T\bigr)W_i\) \hfill \(\triangleright\) Residual outside the current mode-\(i\) space
    \IF{$\|R_i\|_F>\epsilon$}
    \IF{$\texttt{rule}=\texttt{LR-iTT-F}$}
        \STATE $P_i\leftarrow\texttt{orth}(R_i)$
    \ELSE
        \STATE $P_i\leftarrow
        \textsc{SelectReducedBasis}(R_i,\epsilon,\alpha)$
        \hfill $\triangleright$ Algorithm~\ref{alg:reduced_basis}
    \ENDIF
    \STATE $U_i^k\leftarrow
    [\,\widetilde U_i^{k-1}\;\;P_i\,]$,
    $r_i^k\leftarrow r_i^{k-1}+p_i$
    \hfill $\triangleright$ $p_i$ is the number of columns of $P_i$
\ELSE
    \STATE $U_i^k\leftarrow\widetilde U_i^{k-1}$,
    $r_i^k\leftarrow r_i^{k-1}$
\ENDIF
    \STATE \(\mU_i^k\leftarrow
    \texttt{reshape}(U_i^k,[r_{i-1}^k,n_i,r_i^k])\)
    \IF{\(i<d\)}
        \STATE \(\mW_{i+1}\leftarrow
        \texttt{reshape}\bigl((U_i^k)^T W_i V^k_{i+1},
        [r_i^k,n_{i+1},s_{i+1}]\bigr)\) \hfill \(\triangleright\) Contract basis into the next core (uses all columns of \(W_i\))
        \STATE \(W_{i+1}\leftarrow
        \texttt{reshape}(\mW_{i+1},[r_i^k n_{i+1},s_{i+1}])\)
    \ELSE
        \STATE \(W_{d+1}\leftarrow (U_d^k)^T W_d\)
    \ENDIF
\ENDFOR
\STATE Zero-pad \(\mU_{d+1}^{k-1}\) to obtain
\(\widetilde{\mU}_{d+1}^{k-1}\in
\mathbb{R}^{r_d^k\times (k-1)\times 1}\)
\STATE \(\mU_{d+1}^{k}(:,1:k-1,1)\leftarrow
\widetilde{\mU}_{d+1}^{k-1}(:,1:k-1,1)\)
\STATE \(\mU_{d+1}^{k}(:,k,1)\leftarrow W_{d+1}\)
\STATE \(\bpsi^{k}
\leftarrow
\mU^{k}_1\mU^{k}_2\cdots \mU^{k}_{d}\,\mU^{k}_{d+1}\)
\IF{\(\max\limits_{1\le i\le d} r_i^k>r_{\max}\)}
    \STATE \(\bpsi^{k} \leftarrow
    \texttt{TT-rounding}(\bpsi^{k},
    \epsilon_{\rm tt})\)
\ENDIF
\end{algorithmic}
\end{algorithm}

\begin{algorithm}[t]
\footnotesize
\caption{\ssy{\textsc{SelectReducedBasis}: residual basis selection
for LR-iTT-R}}
\label{alg:reduced_basis}
\begin{algorithmic}[1]
\REQUIRE Residual $R\in\mathbb{R}^{m\times s}$ satisfying
$\|R\|_F>\epsilon$; residual tolerance $\epsilon$;
parameter $\alpha\in(0,1]$.
\ENSURE Orthonormal increment $P$ satisfying
$\|(I-PP^T)R\|_F\leq\epsilon$.
\STATE $a_j\leftarrow\|R(:,j)\|_2$ for $j=1,\ldots,s$;
$a_{\max}\leftarrow\max_j a_j$
\STATE $L\leftarrow\{j:a_j<\alpha a_{\max}\}$
\hfill $\triangleright$ Weak columns eligible for removal
\STATE Order the indices in $L$ as
$j_1,\ldots,j_{|L|}$ such that
$a_{j_1}\leq\cdots\leq a_{j_{|L|}}$
\STATE $q\leftarrow
\max\left\{
\ell\in\{0,\ldots,|L|\}:
\sum_{t=1}^{\ell}a_{j_t}^{2}\leq\epsilon^{2}
\right\}$
\hfill $\triangleright$ The empty sum is zero
\STATE $S\leftarrow
\{1,\ldots,s\}\setminus\{j_1,\ldots,j_q\}$
\hfill $\triangleright$ Retain all columns if $q=0$
\STATE $P\leftarrow\texttt{orth}\bigl(R(:,S)\bigr)$
\hfill $\triangleright$ Orthonormal basis for the retained column space
\RETURN $P$
\end{algorithmic}
\end{algorithm}

\subsection{Comparison with existing incremental TT algorithms}
The proposed method is related to deterministic streaming TT algorithms such as
TT-FOA \cite{abed2021adaptive} and TT-ICE \cite{aksoy2024incremental}. Unlike these methods, however, it 
targets streaming low-rank data already represented in TT form. TT-FOA tracks a
fixed-rank TT model of the accumulated data tensor by recursively updating the
cores through a weighted least-squares fit.TT-ICE is closer in spirit to our method. Both methods
can be viewed as tensor extensions of the
projection-then-enrich strategy used in incremental SVD
\cite{brand2006fast}. Specifically, they adaptively expand the TT cores
of an accumulated tensor by appending orthogonal directions
needed to represent incoming tensors. The key distinction
lies in how the incoming data are represented.
TT-FOA uses dense or sampled tensor entries, while TT-ICE assumes access to dense tensor data and forms residuals from their unfoldings. In contrast, the proposed approach takes each incoming tensor directly in TT form and propagates compressed interfaces through a sweep. Thus the projection, residual enrichment, and coefficient update are performed at the TT-core level, without materializing or densely unfolding the incoming tensor. This distinction is essential for data generated by emerging low-rank solvers, where reconstructing dense snapshots would negate the computational advantage of the low-rank representation.

\begin{rem}[Heuristic acceleration]
In addition to TT-ICE, \cite{aksoy2024incremental} introduces
a heuristically accelerated variant, TT-ICE$^*$. Its
heuristics reduce the enrichment cost by skipping basis updates based on the projection
errors of the incoming data and skipping updates
of cores with high occupancy.
Similar strategies could be incorporated into LR-iTT-F
and LR-iTT-R, since snapshot projection errors and core occupancy can be evaluated directly with core-wise operations.
However, heuristic skipping does not automatically preserve
the error bounds established here and would require
additional accuracy checks. 
\end{rem}

\section{\ssy{Theoretical properties}}
In this section, we establish the approximation properties
of LR-iTT-F and LR-iTT-R and analyze their complexities.

\subsection{Approximation properties}

Fix a snapshot $k$. Let $U_i^k=[\widetilde U_i^{k-1}\ P_i]$ be the enriched basis
($U_i^k=\widetilde U_i^{k-1}$ if $P_i=\varnothing$). The Algorithm~\ref{alg:proposed_tt_update} implicitly generates the
intermediate TT tensors
\begin{subequations}
\begin{align}
&\bm{A}^{k}_{0} = \bm{A}^{k},\\
&\bm{A}^{k}_{i} = \mU^{k}_1\cdots \mU^{k}_{i}\,\mathcal{W}_{i+1}\,\mV^k_{i+2}\cdots\mV^k_{d},\quad i=1,\ldots,d-2,  \\
&\bm{A}^{k}_{d-1} = \mU^{k}_1\cdots \mU^{k}_{d-1}\,\mathcal{W}_{d}, \\
&\bm{A}^{k}_{d} = \mU^{k}_1\cdots \mU^{k}_{d}\,\mU^{k}_{d+1}(:,k,1)=\bm{\psi}^{k}(:,:,\cdots,:,k).
\end{align}
\end{subequations}
Define the error increments and the residuals after enrichment by
\begin{equation}
\label{eq:local_error_def}
    \bm E_i:=\bm A_{i-1}^k-\bm A_i^k,
    \qquad
    \widehat R_i:=
    \bigl(I-U_i^k(U_i^k)^T\bigr)W_i,
    \qquad i=1,\ldots,d.
\end{equation}

\begin{lem}[Residual bound after enrichment]
\label{lem:local}
For either LR-iTT-F or LR-iTT-R,
\[
    \widehat R_i=(I-P_iP_i^T)R_i,
    \qquad
    \|\widehat R_i\|_F\leq\epsilon,
\]
where $P_iP_i^T=0$ if no enrichment is performed.
\end{lem}

\begin{myproof}
Since the columns of $P_i$ lie in $\operatorname{range}(R_i)$,
they are orthogonal to the columns of
$\widetilde U_i^{k-1}$, i.e. $(\widetilde{U}_i^{k-1})^TP_i=0$. Using
$U_i^k=[\,\widetilde U_i^{k-1}\;\;P_i\,]$, we obtain
\[
\begin{aligned}
    \widehat R_i
    &=W_i-[\,\widetilde U_i^{k-1}\;\;P_i\,]([\,\widetilde U_i^{k-1}\;\;P_i\,])^TW_i=W_i-\widetilde U_i^{k-1}
      (\widetilde U_i^{k-1})^T W_i-P_iP_i^T W_i\\
    &=R_i-P_iP_i^T R_i=(I-P_iP_i^T)R_i.
\end{aligned}
\]
If $\|R_i\|_F\leq\epsilon$, no enrichment is performed and
$\widehat R_i=R_i$ satisfying $\|\widehat{R}_i\|_F\leq\epsilon$.
Otherwise, enrichment is performed. If LR-iTT-F is applied,  $P_i$ is an orthonormal basis of $\operatorname{range}(R_i)$, yielding $\widehat R_i=0$. If LR-iTT-R is applied, $\|\hat{R}_i\|_F=\|(I-P_iP_i^T)R_i\|_F\leq \epsilon$ by construction
(see Algorithm~\ref{alg:reduced_basis}). Thus,
$\|\widehat R_i\|_F\leq\epsilon$ in all cases.
\end{myproof}

The next lemma shows that the residual bound after enrichment
yields a bound on each error increment.
\begin{lem}[Error increment]
\label{lem:error}
For either LR-iTT-F or LR-iTT-R, define the error increment
at the $i$-th mode by
\[
    \bm E_i:=\bm A_{i-1}^k-\bm A_i^k,
    \qquad i=1,\ldots,d.
\]
Then
\begin{equation}
\label{eq:error}
    \|\bm E_i\|_F\leq\epsilon,
    \qquad i=1,\ldots,d.
\end{equation}
\end{lem}

\begin{myproof}
For $1\leq i<d$, the residual after enrichment is
\[
\widehat R_i
=
\bigl(I-U_i^k(U_i^k)^T\bigr)W_i.
\]
By construction of $\mathcal W_{i+1}$, its right unfolding
is $(U_i^k)^T W_iV_{i+1}^k$. Hence, the corresponding
matricization of $\mathcal U_i^k\mathcal W_{i+1}$ is
$U_i^k(U_i^k)^T W_iV_{i+1}^k$.
Therefore,
\begin{equation}\label{eq:localdiff}
W_iV_{i+1}^k-U_i^k(U_i^k)^T W_iV_{i+1}^k
=
\bigl(I-U_i^k(U_i^k)^T\bigr)W_iV_{i+1}^k
=
\widehat R_iV_{i+1}^k.
\end{equation}

Let $\widehat{\mathcal R}_i
\in\mathbb{R}^{r_{i-1}^k\times n_i\times s_i}$
denote the tensor whose left unfolding is the matrix
$\widehat R_i$. By \eqref{eq:localdiff}, it follows that
\begin{equation}\label{eq:E}
\bm E_i
=
\bm A_{i-1}^k-\bm A_i^k
=
\mathcal U_1^k\cdots\mathcal U_{i-1}^k\,
\widehat{\mathcal R}_i\,
\mathcal V_{i+1}^k\cdots\mathcal V_d^k.
\end{equation}

Since $\mathcal U_1^k,\ldots,\mathcal U_{i-1}^k$
are left-orthogonal and
$\mathcal V_{i+1}^k,\ldots,\mathcal V_d^k$
are right-orthogonal, contraction with these cores
preserves the Frobenius norm. Hence, by
Lemma~\ref{lem:local},
\[
\|\bm E_i\|_F
=
\|\widehat{\mathcal R}_i\|_F
=
\|\widehat R_i\|_F
\leq\epsilon.
\]

For the terminal case $i=d$, by the definition of the
final slice $W_{d+1}$,
\[
W_{d+1}=(U_d^k)^T W_d, \quad\text{hence}, \quad U_d^kW_{d+1}
=
U_d^k(U_d^k)^T W_d.
\]
Therefore, the difference in the left unfoldings of the
$d$-th cores is exactly the residual after enrichment:
\[
W_d-U_d^kW_{d+1}
=
\bigl(I-U_d^k(U_d^k)^T\bigr)W_d
=
\widehat R_d.
\]
Recall that $\bm A_{d-1}^k$ terminates with the core
$\mathcal W_d$, while $\bm A_d^k$ replaces it with
the contraction of $\mathcal U_d^k$ and $W_{d+1}$.
Therefore, the error increment is given by
\[
\bm E_d
=
\bm A_{d-1}^k-\bm A_d^k
=
\mathcal U_1^k\cdots\mathcal U_{d-1}^k\,
\widehat{\mathcal R}_d,
\]
where the residual tensor $\widehat{\mathcal R}_d$
has the left unfolding $\widehat R_d$.
Since $\mathcal U_1^k,\ldots,\mathcal U_{d-1}^k$
are left-orthogonal, we conclude from
Lemma~\ref{lem:local} that
\[
\|\bm E_d\|_F
=
\|\widehat{\mathcal R}_d\|_F
=
\|\widehat R_d\|_F
\leq\epsilon.
\]
This completes the proof.
\end{myproof}

To obtain the error bound for the algorithm, we next establish
the orthogonality of the error increments. We first establish
the following orthogonality propagation property of the
TT decomposition.

\begin{lem}[Orthogonality propagation in TT format]\label{lem:prop}
Let two TT tensors share the same left-orthonormal cores in modes $1,\dots,i-1$. If the left unfoldings of their
$i$-th cores are orthogonal, then the two tensors are orthogonal in the Frobenius inner product.
\end{lem}

\begin{myproof}
See Appendix~\ref{ap:proof}.
\end{myproof}

\begin{lem}[Orthogonality of the error increments]\label{lem:orth}
The increments $\{\bm E_i\}_{i=1}^d$ are mutually orthogonal in the Frobenius inner product:
\begin{equation}\label{eq:orth}
\langle\bm E_i,\bm E_j\rangle_F=0,\qquad i\ne j .
\end{equation}
\end{lem}

\begin{myproof}
Without loss of generality, we consider $1\leq i<j\leq d$.

By \eqref{eq:E}, the first $i-1$ cores of $\bm E_i$ are
$\mathcal U_1^k,\ldots,\mathcal U_{i-1}^k$, which are
left-orthonormal, and its $i$-th core is
$\widehat{\mathcal R}_i$ with left unfolding $\widehat R_i$.

Since $j>i$, $\bm E_j$ has the same first $i-1$ cores,
and its $i$-th core is $\mathcal U_i^k$ with left unfolding
$U_i^k$.

Recall that $\widehat R_i=(I-U_i^k(U_i^k)^T)W_i$ and
$(U_i^k)^T U_i^k=I$, we have
\[
    \widehat R_i^T U_i^k
    =
    \bigl( (I-U_i^k(U_i^k)^T)W_i\bigr)^TU_i^k
    = W_i^T(I-U_i^k(U_i^k)^T)U_i^k=W_i^T(U_i^k-U_i^k)=0.
\]
Lemma~\ref{lem:prop} then implies
$\langle\bm E_i,\bm E_j\rangle_F=0$.
\end{myproof}

We are now ready to establish the following approximation property.

\begin{thm}\label{thm:main}
Let $\bm\psi$ be the accumulated tensor constructed by Algorithm~\ref{alg:proposed_tt_update} from the sequential TT-formatted data
$\{\bm A^k\}_{k=1}^{n_\mu}$, with the enrichment step performed by LR-iTT-F or by LR-iTT-R.
If the final truncation step~\textup{(d)} is not applied, then for each $k=1,\dots,n_\mu$:
\[
\|\bm A^k-\bm\psi(:,:,\cdots,:,k)\|_F\le\sqrt d\,\epsilon;
\]
\end{thm}

\begin{myproof}
Fix $k\in\{1,\dots,n_\mu\}$. When the truncation step~(d) is not performed, the algorithm preserves the previously
constructed accumulated tensor while progressively updating its TT cores. 

In particular,
\[
\bm{\psi}(:,:,\cdots,:,k)=\bm{\psi}^{k}(:,:,\cdots,:,k).
\]
Let \(\{\bm A_i^{k}\}_{i=0}^d\) denote the sequence of intermediate tensors generated by the algorithm, where
\[
\bm A_0^{k}=\bm A^k,
\qquad
\bm A_d^{k}=\bm{\psi}(:,:,\cdots,:,k).
\]
Then
\[
\bm A^k-\bm{\psi}(:,:,\cdots,:,k)
=
\bm A_0^{k}-\bm A_d^{k}
=
\sum_{i=1}^d \bm E_i.
\] 
By Lemma~\ref{lem:error} and Lemma~\ref{lem:orth}, the error increments $\bm E_i$ are bounded $\|\bm E_i\|_F\leq\epsilon$ and mutually orthogonal i.e. $\langle\bm E_i,\bm E_j\rangle_F=0$, then
\[
\|\bm A^k-\bm\psi(:,:,\cdots,:,k)\|_F^2=\sum_{i=1}^d\|\bm E_i\|_F^2=\sum_{i=1}^d\|\widehat{R}_i\|_F^2\le d\,\epsilon^2,
\]
which yields
\[
\|\bm A^k-\bm{\psi}(:,:,\cdots,:,k)\|_F
\le
\sqrt d\,\epsilon.
\]
\end{myproof}

The preceding analysis also yields an exactness property
for LR-iTT-F. When full residual enrichment is performed
at every core, all error increments vanish, and the new
snapshot is represented exactly. This gives the following
corollary.

\begin{cor}[Exactness]\label{cor:exact}
Suppose LR-iTT-F is used. If
$\|R_i\|_F>\epsilon$ for all $i=1,\dots,d$, then the new snapshot is represented exactly:
\begin{equation}\label{eq:exact}
\bm\psi^k(:,:,\cdots,:,k)=\bm A^k .
\end{equation}
\end{cor}
 
\begin{myproof}
If $\|R_i\|_F>\epsilon$, then LR-iTT-F adds the whole span of $R_i$ and $\widehat{R}_i=0$. Therefore, $\|\bm E_i\|_F=\|\widehat{R}_i\|_F=0$ by Lemma \ref{lem:error}.
If this holds for all $i$, then $\bm A^k=\bm A_0^k=\bm A_1^k=\cdots=\bm A_d^k=\bm\psi^k(:,:,\cdots,:,k)$.
\end{myproof}
\subsection{Complexity analysis}
We conclude the section by providing a complexity analysis. 
Assume that the accumulated TT tensor \(\psi^{(k-1)}\) has maximum rank \(r\), and the \(k\)-th incoming tensor \(A^{(k)}\) has maximum rank \(s\). Assume also that the mode sizes satisfy \(n_i\le n\). During one incremental update, the enrichment at each mode can append at most \(s\) new vectors, and hence the intermediate updated TT ranks are bounded by \(
\hat{r} := r+s\).

The storage of the accumulated TT tensor \(\psi^{(k)}\) is bounded by \(O(dn\hat{r}^2+k\hat{r})\), where the first term corresponds to the physical TT cores and the second term
comes from the final sample core of size \(r_d^{(k)}\times k\times 1\).  Thus, the method avoids any dense tensor storage and scales linearly with the tensor dimension \(d\) when the TT ranks remain moderate.

We next estimate the cost of one streaming update without intermediate truncation. At mode \(i\), since
\(
\widetilde U_i^{(k-1)}\in
\mathbb{R}^{r_{i-1}^{(k)}n_i\times r_i^{(k-1)}},
\)
and
\(
W_i\in \mathbb{R}^{r_{i-1}^{(k)}n_i\times s_i},
\)
and \(r_{i-1}^{(k)}\le \hat{r}\), \(r_i^{(k-1)}\le r\), and \(s_i\le s\), the residual projection
\(
R_i=(I-\widetilde U_i^{(k-1)}(\widetilde U_i^{(k-1)})^T)W_i
\)
costs \(O(n_i\hat{r}r s)\), which is bounded by \(O(n\hat{r}^2 s)\). If enrichment is triggered, computing \(P_i=\texttt{orth}(R_i)\) by a QR factorization costs \(O(n_i\hat{r} s^2)\). The propagation step
\(
W_{i+1}
=
\texttt{reshape}\big((U_i^{(k)})^T W_i V_{i+1}^{(k)}\big)\)
costs \(O(n_i\hat{r}^2s+n_{i+1}\hat{r} s^2)\), which is bounded by \(O(n\hat{r}^2s)\), as \(\hat{r}\ge s\). Therefore, the per-mode cost is bounded by \(O(n\hat{r}^2s)\), and summing over all modes gives
$O\left(dn\hat r^2s\right)$.

For LR-iTT-R, the column norms of $R_i$ are obtained from the residual already formed, at cost $O(n_i\hat rs)$; sorting the weak columns costs $O(s\log s)$; and the QR factorization is applied to $R_i(:,S_i)$ with $|S_i|\le s$ columns, hence its cost is bounded by $O(n_i\hat rs^2)$ as well. Moreover, $p_i\le|S_i|\le s$, so the rank bound $\hat r=r+s$ and the storage estimate remain valid. Consequently the per-mode cost $O(n\hat r^2s)$ and the total cost $O(dn\hat r^2s)$ are unchanged.

If intermediate \texttt{TT-rounding} is triggered, the additional cost is
$O\left(dn \hat r^3+k \hat r^2\right)$,
where the first term corresponds to rounding the physical TT cores and the second term comes from the final core.

Thus, for moderate TT ranks, Algorithm 3.1 scales linearly with the tensor dimension \(d\) in both storage and computational cost.

\section{Model order reduction with the incremental TT reduced basis\label{sec:tt_rom}}

As established, the proposed incremental algorithm produces a
compressed accumulated TT representation of streaming
TT-formatted data. We now show how this representation can
be used for model order reduction of linear kinetic equations.
The accumulated TT tensor plays a role analogous to a proper
orthogonal decomposition (POD) representation of training data:
it simultaneously encodes a compressed reduced basis and the
corresponding coefficients associated with the training
parameters.

In this section, we first review classical POD for full-order
data and explain its relation to the accumulated TT
representation. We then show how the resulting compressed
TT reduced basis can be used to construct ROMs directly from
low-rank solution data through core-level operations, in both
non-intrusive and intrusive settings.

\subsection{Classical proper orthogonal decomposition (POD)}
Consider a parameterized discrete full-order model
\begin{equation}
\label{eq:classical_fom}
    A(\mu)u(\mu)=f(\mu), \qquad u(\mu)\in \mathbb{R}^{N},
\end{equation}
where \(\mu\in\mathcal{P}\) denotes the parameter such as configurations of material properties or boundary conditions. POD method follows an offline-online decomposition framework. In the offline stage, a reduced space
is constructed by extracting low-rank structures across parameters from the training data. In the online stage, the prediction for new parameter can be made through an interpolation inside or a projection onto the reduced space constructed offline.

In the offline stage, given  a prescribed  error tolerance $\epsilon$ and the snapshot matrix
\begin{equation}
\label{eq:classical_snapshot_matrix}
    S =
    \big[ u(\mu_1),\,u(\mu_2),\,\ldots,\,u(\mu_{n_\mu}) \big]
    \in \mathbb{R}^{N\times n_\mu},
\end{equation}
the POD method constructs the reduced basis by the singular value decomposition:
\begin{equation}
\label{eq:classical_pod_svd}
    S \approx U_{r}\Sigma_{r} V_{r}^T\quad\textrm{with}\quad    \|S-U_{r}\Sigma_{r}V_{r}^T\|_F
    \le \epsilon.
\end{equation}
Here, the reduced dimension $r$ is selected as
\begin{equation}
\label{eq:classical_pod_rank}
    r
    =
    \min\left\{
    m:\sum_{j=m+1}^{\min(N,n_\mu)}\sigma_j^2
    \le \epsilon^2
    \right\},
\end{equation}
where \(\{\sigma_j\}\) are the singular values of \(S\). 
The columns of \(U_{r}\in\mathbb{R}^{N\times r}\) form the
reduced basis, while the columns of
\(\Sigma_{r}V_{r}^T\in\mathbb{R}^{r\times n_{\mu}}\) provide the reduced coefficients for
the training parameters. In the online stage, for a new parameter \(\mu\), a POD-based ROM
seeks an approximation of the form
\begin{equation}
\label{eq:classical_pod_expansion}
    u(\mu)\approx U_{r}c(\mu),
    \qquad c(\mu)\in\mathbb{R}^{r}.
\end{equation}
The reduced coefficients $c(\mu)$ may be obtained non-intrusively by interpolations in the
parameter space, or intrusively by projecting the governing equations onto the
reduced space spanned by $U_r$. These two approaches will be detailed after introducing the
corresponding TT basis and coefficient core.

\subsection{TT reduced basis and coefficient core}

In the high-dimensional tensor setting, the full-order state is viewed as a
tensor
\[
    u(\mu)\in\mathbb{R}^{n_1\times n_2\times\cdots\times n_d},
    \qquad N=n_1n_2\cdots n_d.
\]
Classical POD first vectorizes each snapshot and forms the 
snapshot matrix
\[
    S=\big[\operatorname{vec}(u(\mu_1)),\ldots,
    \operatorname{vec}(u(\mu_{n_\mu}))\big]
    \in\mathbb{R}^{(n_1n_2\cdots n_d)\times n_\mu}.
\]
Equivalently, \(S\) is the unfolding that separates all physical DOFs from the training sample index.

For simplicity, we treat the training samples as a single additional mode indexed by \(k\); that is, the parameter sampling space is not further
tensorized in the following construction. The
proposed method constructs the accumulated snapshot tensor
\begin{equation}
\label{eq:tt_snapshot_tensor_rom}
    \bm{\psi}(i_1,\ldots,i_d,k)
    \approx
    u(\mu_k)(i_1,\ldots,i_d),
    \qquad k=1,\ldots,n_\mu,
\end{equation}
directly in TT form,
\begin{equation}
\label{eq:tt_snapshot_tensor_decomp_rom}
    \bm{\psi}
    =
    \mU_1\mU_2\cdots \mU_d\mU_{d+1}.
\end{equation}
Unfolding \(\bm{\psi}\) across the physical-sample interface
\((i_1,\ldots,i_d)\mid k\) gives
\begin{equation}
\label{eq:tt_physical_sample_unfolding}
    S_{\rm TT}
    =
    \bm{\psi}_{\langle d\rangle}
    \in
    \mathbb{R}^{(n_1n_2\cdots n_d)\times n_\mu}.
\end{equation}
This unfolding gives a POD-like factorization of the classical snapshot matrix:
\begin{equation}
\label{eq:tt_unfolding_factorization}
    S\approx S_{\rm TT}=\Phi C,
    \quad
    \Phi\in\mathbb{R}^{(n_1n_2\cdots n_d)\times r_d},
    \quad
    C\in\mathbb{R}^{r_d\times n_\mu}.
\end{equation}
Here, \(\Phi\) denotes the unfolding of compressed reduced basis represented implicitly by
the first \(d\) TT cores.  More precisely, its tensorization satisfies
\begin{equation}
\label{eq:tt_basis_rom}
    \boldsymbol{\Phi}(i_1,\ldots,i_d,a)
    =
    \mU_1(1,i_1,:)\,
    \mU_2(:,i_2,:)\cdots
    \mU_d(:,i_d,a),
    \qquad a=1,\ldots,r_d .
\end{equation}
Since the physical cores produced by the incremental construction
are left-orthogonal, this basis is orthonormal, i.e.
\(\Phi^T\Phi=I_{r_d}\). 
The coefficient matrix \(C\) is
encoded by the final sample core,
\begin{equation}
\label{eq:tt_coeff_core_rom}
    C(a,k)=c_a(\mu_k)=\mU_{d+1}(a,k,1),
    \qquad a=1,\ldots,r_d,\quad k=1,\ldots,n_\mu.
\end{equation}
Crucially, the unfolding is only conceptual: in practice, \(\Phi\) is never
explicitly constructed as a dense matrix, but is stored and manipulated through
its TT cores.

In summary, the first $d$ physical cores defines a  TT analogue of the POD basis, while
\(C\) plays the role of the POD coefficient matrix \(\Sigma_r V_r^T\). In contrast to classical POD, neither the full snapshots, the full snapshot
matrix, nor the dense reduced basis matrix is formed.

\begin{rem}
The use of a single sample index \(k\) is only for notational simplicity. If
the training parameters are sampled on a tensor-product parameter grid, the
sample mode may itself be decomposed into several parameter modes. In that
case, the accumulated snapshot tensor has the form
\[
    \bm{\psi}(i_1,\ldots,i_d,\ell_1,\ldots,\ell_p)
    \approx
    u(\mu_{\ell_1,\ldots,\ell_p})(i_1,\ldots,i_d),
\]
and the same construction applies by taking the interface between the physical
modes \((i_1,\ldots,i_d)\) and the parameter/sample modes
\((\ell_1,\ldots,\ell_p)\). The reduced basis is determined by the TT rank
across this interface, while the remaining parameter cores encode the training
coefficients in compressed form.
\end{rem}

\subsection{Non-intrusive coefficient interpolation\label{sec:non-intrusive}}
In the non-intrusive approach, solutions at new parameter values are predicted
by interpolating or regressing the coefficient core, without requiring access to
the governing operators. From
\eqref{eq:tt_coeff_core_rom}, the training data provide reduced coefficient
vectors
\begin{equation}
\label{eq:tt_training_coefficients}
    c(\mu_k)=\mU_{d+1}(:,k,1)\in\mathbb{R}^{r_d},
    \qquad k=1,\ldots,n_\mu.
\end{equation}
A non-intrusive surrogate is obtained by constructing an interpolant or
regression model for $\mu\mapsto c(\mu)$.
For example, one may use radial basis functions, polynomial or sparse-grid
interpolations, Gaussian regressions, or neural-network regressions,
depending on the dimension and structure of the parameter space.

For a new parameter \(\mu_\ast\), the interpolated or regressed coefficient
\(c(\mu_\ast)\) defines the reduced approximation
\begin{equation}
\label{eq:tt_nonintrusive_rom}
    u_r(\mu_\ast)
    =
    \sum_{a=1}^{r_d}\Phi_a\,c_a(\mu_\ast).
\end{equation}
The reconstruction can also remain in TT form. Specifically, the coefficient
\(c(\mu_\ast)\) is contracted into the last rank index of the block TT
basis \eqref{eq:tt_basis_rom}.

\subsection{Intrusive method based on Galerkin projection}

In the intrusive approach, solutions to a parametric linear system $A(\mu)u(\mu)=f(\mu)$ at new parameter values are obtained by
projecting the governing equations onto the TT reduced basis, which requires
access to the full-order operators. We assume the operator and the right-hand side admit affine decompositions
\begin{equation}
\label{eq:tt_affine_operator_rhs}
    A(\mu)=\sum_{q=1}^{Q}\theta_q(\mu)A_q,
    \qquad
    f(\mu)=\sum_{s=1}^{S}\eta_s(\mu)f_s,
\end{equation}
where each \(A_q\) has tensor-product structure
\begin{equation}
\label{eq:tt_tensor_product_operator}
    A_q
    =
    A_{q,1}\otimes A_{q,2}\otimes\cdots\otimes A_{q,d},
    \qquad A_{q,i}\in\mathbb{R}^{n_i\times n_i},
\end{equation}
and each \(f_s\) is given in separable or TT form.  Under these assumptions,
the reduced operator and right-hand side can be assembled by TT contractions
without forming the full matrices or vectors.

Denote the unfolding of the TT reduced basis as $\Phi$. For a new parameter, we prediction the solution by solving the reduced Galerkin system:
\begin{subequations}
\label{eq:tt_galerkin_rom}
\begin{align}
\label{eq:tt_reduced_system}
   & u(\mu)\approx\Phi c(\mu),
    \qquad
    A_r(\mu)c(\mu)=f_r(\mu),
\\
\label{eq:tt_reduced_operator_chain}
    &A_r(\mu)
    =
    \Phi^T A(\mu)\Phi
    =
    \sum_{q=1}^{Q}\theta_q(\mu)A_{r,q},
    \qquad
    A_{r,q}=\Phi^T A_q\Phi,\\
\label{eq:tt_reduced_rhs_chain}
   & f_r(\mu)
    =
    \Phi^T f(\mu)
    =
    \sum_{s=1}^{S}\eta_s(\mu)f_{r,s},
    \qquad
    f_{r,s}=\Phi^T f_s.
\end{align}
\end{subequations}
We now describe how each reduced matrix \(A_{r,q}\) is computed without
assembling either the full operator \(A_q\) or the full basis matrix \(\Phi\),
using a left-to-right TT contraction.

\begin{enumerate}
\item \textbf{Local matricization of the TT basis core.}
 For each physical mode
\(i=1,\ldots,d\), let \(U_i\) denote the left unfolding of the
\(i\)-th TT basis core,
\[
    U_i
    =
    \texttt{reshape}\!\left(\mU_i,\,[r_{i-1}n_i,\;r_i]\right),
    \qquad i=1,\ldots,d.
\]
Here \(r_0=1\), and \(r_d\) is the reduced basis dimension. The columns of
\(U_i\) are orthonormal since the basis cores are left-orthogonal.

\item \textbf{Left-to-right partial contraction sweep.}
For a fixed affine component \(A_q\), initialize the partial contraction matrix
as $E_0^{(q)}=[1].$
The matrix \(E_{i-1}^{(q)}\in\mathbb{R}^{r_{i-1}\times r_{i-1}}\)
contains the partial Galerkin contraction over modes \(1,\ldots,i-1\). At mode
\(i\), first combine this partial contraction matrix with the local physical
operator \(A_{q,i}\):
\begin{equation}
\label{eq:tt_galerkin_local_operator}
    B_i^{(q)}
    =
    E_{i-1}^{(q)}\otimes A_{q,i}
    \in
    \mathbb{R}^{(r_{i-1}n_i)\times(r_{i-1}n_i)}.
\end{equation}
Then project through the left unfolding of the current core, $U_i\in\mathbb{R}^{r_{i-1}n_i\times r_i}$, advancing the 
interface rank from \(r_{i-1}\) to \(r_i\):
\begin{equation}
\label{eq:tt_galerkin_environment_update}
    E_i^{(q)}
    =
    U_i^T
    B_i^{(q)}
    U_i\in\mathbb{R}^{r_{i}\times r_i}.
\end{equation}
	After the last physical mode, all physical indices have been contracted and
	the reduced matrix is formed as
\begin{equation}
\label{eq:tt_galerkin_reduced_matrix}
	    A_{r,q}=E_d^{(q)}\in\mathbb{R}^{r_d\times r_d}.
		\end{equation}
		\end{enumerate}

For affine problems, the reduced components \(A_{r,q}\) and \(f_{r,s}\) are
precomputed and stored offline.  For a new parameter \(\mu\), the online stage
assembles the reduced system by linear combination of these components and then
solves \eqref{eq:tt_galerkin_rom}.

\begin{rem}
The affine tensor-structured assumptions above are crucial for an efficient
offline-online decomposition. In the full-rank setting, when the operator
\(A(\mu)\) or the right-hand side \(f(\mu)\) is non-affine, one may first
construct an approximate affine surrogate using the empirical interpolation
method (EIM) \cite{barrault2004empirical} or discrete EIM (DEIM)
 \cite{chaturantabut2010nonlinear}. These techniques can, in
principle, be extended to low-rank formats, in the same spirit as the extension
of matrix DEIM-CUR \cite{sorensen2016deim} to tensor-valued problems
\cite{dektor2024collocation}. The focus of this paper is the incremental basis
construction and its application to ROM; such hyper-reduction extensions are
therefore deferred to future work.
\end{rem}

\section{Numerical experiments}
\label{sec:numerics}

This section evaluates the proposed incremental TT construction on
parametric RTE data that are already available in low-rank form. The experiments are designed to assess two
aspects of the proposed methodology: (i) streaming TT data compression and (ii) its application in the intrusive projection-based and non-intrusive interpolation-based ROMs.

For the compression study, we compare our approach with TT-ICE
\cite{aksoy2024incremental} using three diagnostics: reconstruction error, TT
ranks of the accumulated solution tensor, and wall time. In particular, for an accumulated tensor \(\bpsi\) produced by either method from the streaming data \(\{\bm A^k\}_{k=1}^{n_\mu}\), we define the relative reconstruction error of the $k$-th snapshot by
\[
 e_k= \frac{\|\bm A^k-\bpsi(:,\ldots,:,k)\|_F}{\|\bm A^k\|_F},
\qquad k=1,\ldots,n_\mu,
\]
and report the mean reconstruction error
\[\frac{1}{n_\mu}\sum_{k=1}^{n_\mu} e_k.\]
The proposed method receives each snapshot as a TT tensor and updates the accumulated tensor
directly through operations on the TT cores. In contrast, TT-ICE operates on the dense representation, requiring each incoming snapshot to be materialized as a full tensor. Both methods process the same training snapshots in the same order, while using their respective input representations, and employ the same
intermediate TT-rounding strategy. We note that the original TT-ICE algorithm in \cite{aksoy2024incremental} does not perform reorthogonalization or intermediate truncation during the iterative updates, primarily to avoid the associated computational overhead. In our
experiments, however, we observe that moderate intermediate truncation is
beneficial for all three methods, since it removes redundant information and
prevents excessive growth of the TT ranks of the accumulated tensor. Although
such truncation introduces additional cost, it can substantially reduce memory
usage and the cost of subsequent operations. Therefore, in the comparisons
below, we employ the truncation strategy for all three methods as described in Algorithm~\ref{alg:proposed_tt_update}: the TT ranks of
the accumulated tensor are monitored during the updates, and \texttt{TT-rounding} is performed only when the maximum rank exceeds a prescribed threshold. The truncation tolerance is chosen according to the target accuracy in each experiment.

For the ROM study, we incrementally construct a TT reduced basis with the proposed method to build both intrusive and non-intrusive ROMs, following the formulations presented in Section~\ref{sec:tt_rom}. We demonstrate the accuracy and online efficiency of the ROMs, as well as the
offline computational savings obtained by using low-rank rather than full-rank
training data. All numerical experiments were performed in MATLAB on a MacBook Air equipped with an Apple M1 chip and 16 GB of RAM.

The remainder of this section is organized as follows. Section~\ref{sec:benchmark}
introduces the benchmark examples. Section~\ref{sec:compression}
compares the performance of the proposed method with TT-ICE in data compression. Section~\ref{sec:numerical-rom} investigates the use of the proposed incremental TT compression in ROMs. 

\subsection{Setup of benchmark examples\label{sec:benchmark}}
We consider three 2D2V benchmark examples of the steady-state
parametric liner kinet equation~\eqref{eq:rte}, with spatial location
$\bx=(x,y)\in\Gamma_\bx\subset\mathbb{R}^2$ and angular
direction $\BOmega\in\mathbb{S}^2$. The parameter dependence
enters through the scattering and absorption cross sections,
while the source is independent of the parameters.
Homogeneous inflow boundary conditions are imposed in all
examples. The domains, material configurations, sources,
and training parameters are specified below.

We generate the training snapshots using the sweep-based
low-rank source iteration with diffusion synthetic
acceleration in \cite{guo2026highly}. The solver uses an
upwind discontinuous Galerkin spatial discretization with
first-order polynomial approximation and a discrete
ordinates angular discretization based on the
Chebyshev--Legendre quadrature
$\mathrm{CL}(N_\theta,N_v)$
\cite{lewis1983computational}. It directly computes the
solution in low-rank form. We use the default convergence
criterion specified in the numerical section of
\cite{guo2026highly}.

\textbf{Homogeneous medium.}
The homogeneous example is posed on \(\Gamma_\bx=[-1,1]^2\). There is no absorption, and the parametric scattering
cross section is $\sigma_s(x,y;\mu_s)=\mu_s$ with $\mu_s\in[85,105]$. A Gaussian source is imposed:
$G(x,y)=\exp\!\left(-100(x^2+y^2)\right).$
 
We use a $64\times64$ uniform rectangular mesh. The angular domain is discretized using the \(\mathrm{CL}(32,16)\) rule,
giving TT mode sizes \((128,128,16,32)\). The \(41\) training parameters are uniformly sampled as
\(\mu_s\in\{85,85.5,\ldots,105\}\).   

\begin{figure}[t]
\centering
\begin{tabular}{@{}cc@{}}
\includegraphics[width=0.44\textwidth]{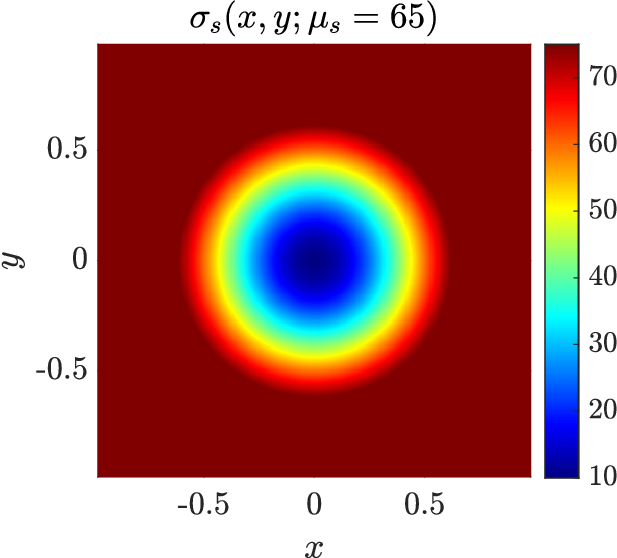} &
\includegraphics[width=0.44\textwidth]{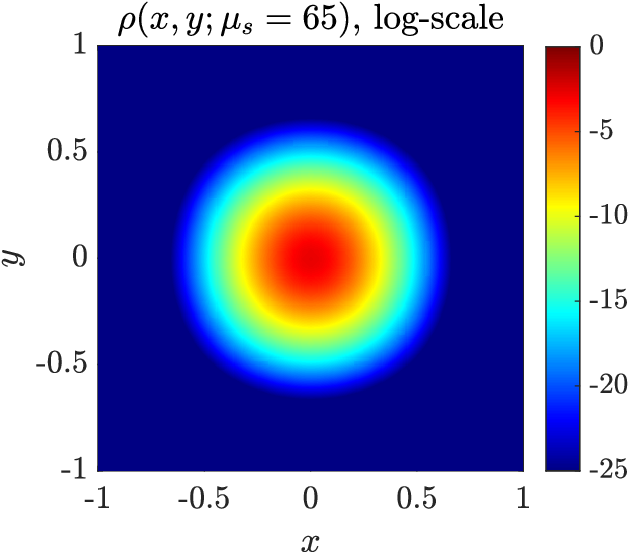} \\
\small (a) Changing-scattering material, \(\mu_s=65\) &
\small (b) Changing-scattering scalar flux, \(\mu_s=65\) \\[0.8em]
\includegraphics[width=0.38\textwidth]{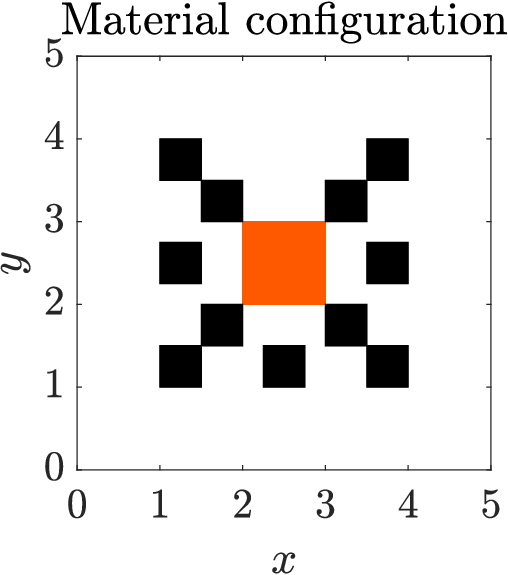} &
\includegraphics[width=0.44\textwidth]{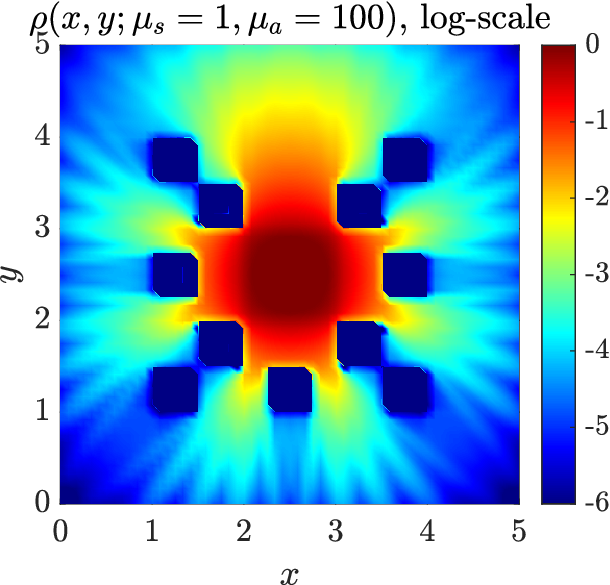} \\
\small (c) Lattice material/source configuration &
\small (d) Lattice scalar flux, \((\mu_s,\mu_a)=(1,100)\)
\end{tabular}
\caption{Representative material configurations and corresponding scalar fluxes for the
changing-scattering and lattice benchmarks.  Panels (a) and (c) show the material configuration, while panels (b) and (d) show
\(\rho(x,y)\) under the logarithmic scale, computed from the stored low-rank
training snapshots.}
\label{fig:rte-benchmark-visuals}
\end{figure}

\textbf{Changing-scattering medium.}
The changing-scattering example is also posed on \(\Gamma_\bx=[-1,1]^2\), with
the Gaussian source $G(x,y)=\exp\!\left(-100(x^2+y^2)\right).$
The parametric scattering cross section is defined as 
\begin{equation}
\label{eq:changing-scattering-benchmark}
    \sigma_s(x,y;\mu_s)=\begin{cases}
    \mu_s r^2(r^2-2)^2+10, & r\le 1,\\
    \mu_s+10, & r>1.
    \end{cases},
    \qquad
    \sigma_a(x,y;\mu_s)=10,
\end{equation}
where $\mu_s\in[40,90]$ and $r=\sqrt{x^2+y^2}$. We present the reference solution and the material configuration corresponding to $\mu_s=65$ in the top row of Figure~\ref{fig:rte-benchmark-visuals}.

We use a \(50\times 50\) spatial mesh and
\(\mathrm{CL}(32,16)\) angular discretization, resulting in mode sizes
\((100,100,16,32)\). The training set contains \(51\) uniformly spaced samples
\(\mu_s\in\{40,41,\ldots,90\}\). 

\textbf{Lattice medium.}
The lattice example follows the two-material configuration on \(\Gamma_\bx=[0,5]^2\) shown in panel (c) of Figure~\ref{fig:rte-benchmark-visuals}.
The black region is purely absorbing, the other regions are purely scattering, and a constant source is imposed in the centered orange region.
The parametric scattering and absorption cross sections are defined as piecewise constants:
\begin{equation}
\label{eq:lattice-cross-sections}
\begin{aligned}
    (\sigma_a(x,y;\mu_a),\sigma_s(x,y;\mu_s))
    =
    \begin{cases}
    (\mu_a,0), & (x,y)\in \textrm{black regions},\\
    (0,\mu_s), &  \textrm{otherwise}.
    \end{cases}.
\end{aligned}
\end{equation}
Here, we consider $(\mu_a,\mu_s)\in[95,105]\times[0.5,1.5]$. We present the reference solution corresponding to $(\mu_a,\mu_s)=(100,1)$ in the bottom row of Figure~\ref{fig:rte-benchmark-visuals}.

We use a \(60\times 60\) uniform spatial mesh and \(\mathrm{CL}(24,12)\) for the angular discretization, resulting in TT mode sizes \((120,120,12,24)\).
In addition, we generate training data using an $11\times 11$ uniform sampling in the parameter domain $[95,105]\times[0.5,1.5]$.

For convenience, we summarize the benchmark configurations in Table~\ref{tab:rte-benchmark-configs}. 
\begin{table}[t]
\caption{Benchmark configurations for the radiative-transfer examples.}
\label{tab:rte-benchmark-configs}
\centering
\small
\resizebox{\textwidth}{!}{%
\begin{tabular}{cccccc}
\toprule
Example & Domain & \(N_x\times N_y\) & \(N_{v_z}\times N_\theta\) &
TT mode sizes & Training parameters \\
\midrule
Homogeneous & \([-1,1]^2\) & \(64\times 64\) & \(16\times 32\) &
\((128,128,16,32)\) & \(\mu_s\in\{85,85.5,\ldots,105\}\) \\
Changing scattering & \([-1,1]^2\) & \(50\times 50\) & \(16\times 32\) &
\((100,100,16,32)\) & \(\mu_s\in\{40,41,\ldots,90\}\) \\
Lattice & \([0,5]^2\) & \(60\times 60\) & \(12\times 24\) &
\((120,120,12,24)\) &
\(\mu_s\in\{0.5,0.6,\ldots,1.5\}\), \(\mu_a\in\{95,96,\ldots,105\}\) \\
\bottomrule
\end{tabular}
}
\end{table}

\subsection{Data compression}
\label{sec:compression}
We use the homogeneous and changing-scattering examples with the setup from
Table~\ref{tab:rte-benchmark-configs} to evaluate the proposed method as a
streaming TT compression procedure.  

\subsubsection{Homogeneous medium}

We first consider the homogeneous example with the configuration given in Table~\ref{tab:rte-benchmark-configs}. The low-rank solver \cite{guo2026highly} generates a collection of 41 solution snapshots in low-rank matrix-factorized form. Then we construct three TT datasets by converting each snapshot to TT format using conversion tolerances \(10^{-4}\), \(10^{-5}\), and \(10^{-6}\), respectively. After conversion, each TT snapshot has mode size
\(n=(128,128,16,32)\)
and TT rank structure \((1,r_1,r_2,r_3,1)\). The corresponding mean internal TT ranks \((r_1,r_2,r_3)\) of the three TT datasets are
$$(13.20,14.98,11.37),\qquad
 (21.24,26.59,17.49),\qquad
 (28.68,42.68,24.05),$$
respectively.

\textbf{Reconstruction accuracy and wall time.}
\pzc{For each of the three tolerance settings, all methods use the same nominal
relative tolerance $\epsilon$, TT-rounding tolerance
$\epsilon_{\mathrm{tt}}$, and rank trigger $r_{\max}$ listed in
Table~\ref{tab:two-d-two-v-results}. We omit the final TT-rounding step,
so the reported ranks reflect the streaming updates, including any
intermediate rounding operations.}

\ssy{Table~\ref{tab:two-d-two-v-results} summarizes the mean reconstruction
errors, TT ranks of the accumulated tensor, and wall times. Both proposed
methods achieve lower mean reconstruction errors than TT-ICE in all
three settings. For $\epsilon=10^{-5}$,
$\epsilon_{\mathrm{tt}}=10^{-6}$, and $r_{\max}=30$, LR-iTT-R
achieves a mean error of $7.092\times10^{-6}$, compared with
$1.566\times10^{-5}$ for TT-ICE. At the most stringent setting,
the corresponding errors are $7.110\times10^{-7}$ and
$2.943\times10^{-6}$, respectively. Both proposed methods also
require substantially less wall time than TT-ICE. In the order
listed in the table, LR-iTT-F achieves speedups of
$2.44\times10^{2}$, $1.00\times10^{2}$, and $1.90\times10^{1}$,
while LR-iTT-R achieves speedups of $2.74\times10^{2}$,
$1.19\times10^{2}$, and $1.98\times10^{1}$. The proposed methods generally produce slightly larger ranks between physical
modes than TT-ICE, but  smaller ranks at the physical-sample interface,
yielding lower-dimensional ROMs.}


Furthermore, the speedup decreases as the tolerances are tightened, since the TT ranks of the accumulated tensor \(\bpsi\) increase, leading to higher costs for the projection. In contrast, the runtime of TT-ICE is dominated by residual computations and the associated SVDs. Because the dimensions of these operations are determined primarily by the incoming snapshot size, they are less sensitive to the rank growth of the accumulated tensor. Consequently, the computational cost of the proposed methods increase more rapidly as the accuracy requirements become more stringent, reducing the observed speedup while still maintaining a significant performance advantage.

\begin{figure}[t]
\centering
\includegraphics[width=\textwidth]{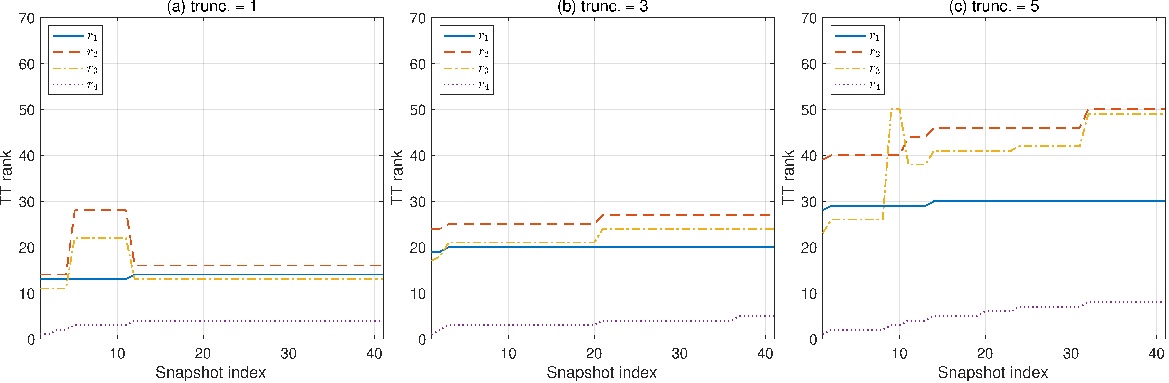}
\par\vspace{4pt}
\includegraphics[width=\textwidth]{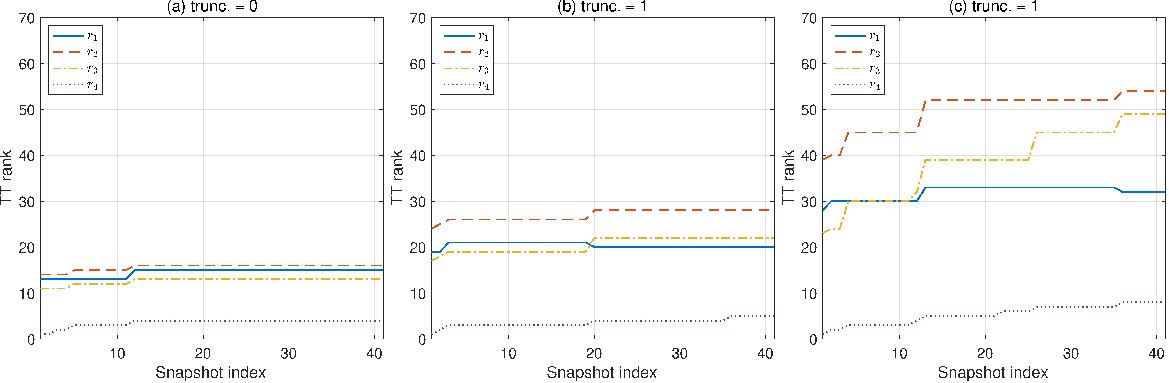}
\par\vspace{4pt}
\includegraphics[width=\textwidth]{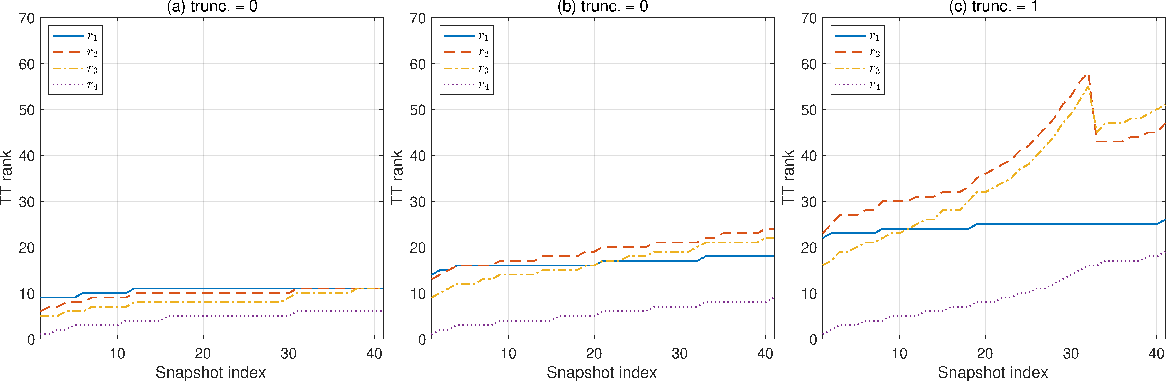}
\caption{\ssy{Evolution of the four internal TT ranks on the homogeneous
2D2V training set under the three parameter settings in
Table~\ref{tab:two-d-two-v-results}. The top, middle, and bottom rows
correspond to LR-iTT-F, LR-iTT-R, and TT-ICE, respectively.
In each row, the panels correspond to:
(a) $\epsilon=10^{-4}$, $\epsilon_{\rm tt}=10^{-5}$, $r_{\max}=30$;
(b) $\epsilon=10^{-5}$, $\epsilon_{\rm tt}=10^{-6}$, $r_{\max}=30$;
(c) $\epsilon=10^{-6}$, $\epsilon_{\rm tt}=10^{-7}$, $r_{\max}=60$.}}
\label{fig:homogeneous-rank-evolution}
\end{figure}

\textbf{TT-rank evolution.}
\ssy{In Figure \ref{fig:homogeneous-rank-evolution}, we plot the evolution of the four internal TT ranks of the accumulated tensor \(\bpsi\), together with the number of intermediate rounding operations for LR-iTT-F, LR-iTT-R, and TT-ICE. Generally, TT-ICE has smaller ranks between the physical modes, whereas its rank at the physical-sample interface is larger than that of either proposed method. Across the three settings, LR-iTT-F performs 1, 3, and 5 intermediate rounding operations, LR-iTT-R performs 0, 1, and 1, and TT-ICE performs 0, 0, and 1, respectively. Thus, intermediate rounding contributes to the observed rank histories, particularly for LR-iTT-F at tighter tolerances.}

\begin{figure}[t]
\centering
\includegraphics[width=\textwidth]{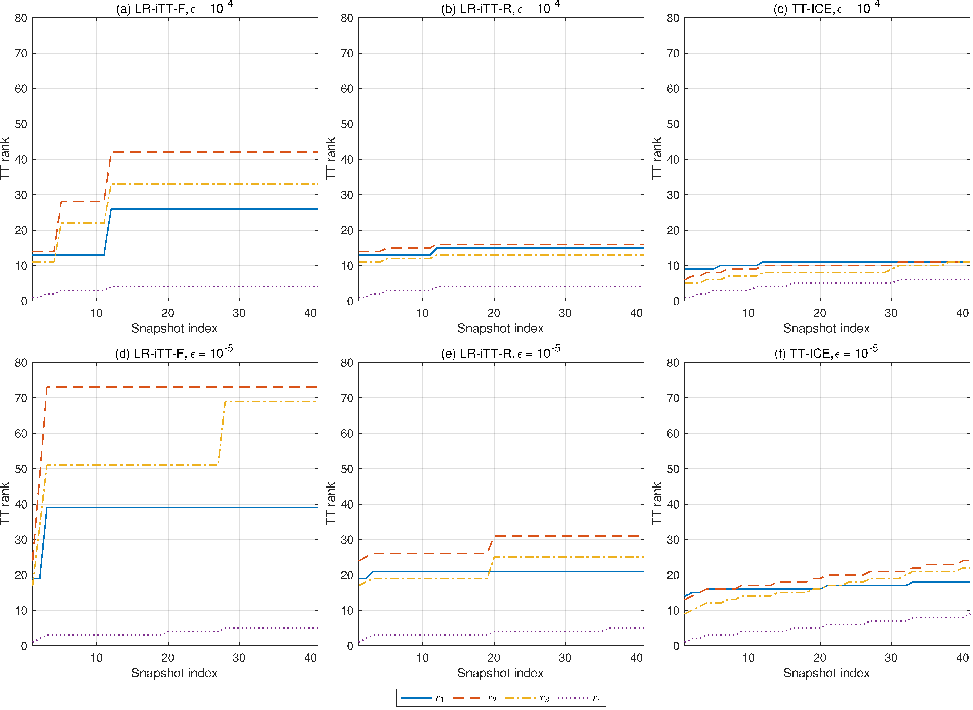}
\caption{\ssy{Comparison on the homogeneous 2D2V training set
without intermediate truncation. Panels (a), (b) and (c) use
\(\epsilon=10^{-4}\), while panels (d), (e) and (f) use \(\epsilon=10^{-5}\). Numerical results are report in Table~\ref{tab:two-d-two-v-without-truncation}.}}
\label{fig:homogeneous-no-truncation-rank-comparison}
\end{figure}

\textbf{Effect of  TT rounding.}
\ssy{To assess the effect of intermediate TT-rounding, we rerun the compression test with $(\epsilon,\epsilon_{\mathrm{tt}})=(10^{-4},10^{-5})$ and $(10^{-5},10^{-6})$ without prescribing the maximum rank. In other words, intermediate and final TT-rounding are disabled. Results are shown in Figure~\ref{fig:homogeneous-no-truncation-rank-comparison} and Table~\ref{tab:two-d-two-v-without-truncation}. In these runs, all internal TT ranks are non-decreasing. Relative to Table~\ref{tab:two-d-two-v-results}, the final ranks of LR-iTT-F increase from $(14,16,13,4)$ to $(26,42,33,4)$ at $\epsilon=10^{-4}$, and from $(20,27,24,5)$ to $(39,73,69,5)$ at $\epsilon=10^{-5}$. The impact of removing rounding on LR-iTT-R is much milder: its ranks remain $(15,16,13,4)$ at $\epsilon=10^{-4}$ and increase modestly from $(20,28,22,5)$ to $(21,31,25,5)$ at $\epsilon=10^{-5}$. The TT-ICE ranks are unchanged in both cases because no intermediate rounding was triggered in the corresponding runs in Table~\ref{tab:two-d-two-v-results}. Without intermediate rounding, both proposed methods still achieve
lower mean reconstruction errors and shorter compression times than
TT-ICE. At $\epsilon=10^{-4}$ and $10^{-5}$, LR-iTT-R obtains errors of $3.903\times10^{-5}$ and $6.667\times10^{-6}$, compared with $1.509\times10^{-4}$ and $1.566\times10^{-5}$ for TT-ICE. LR-iTT-F achieves slightly lower errors of $3.283\times10^{-5}$ and $5.445\times10^{-6}$, respectively. LR-iTT-R nevertheless has the lowest average wall time in both cases. It is $3.18$ and $28.56$ times faster than LR-iTT-F, and about $170$ and $106$ times faster than TT-ICE. Without rounding, the ranks between physical mode  is substantially larger than TT-ICE with LR-iTT-F,  but only slightly larger than TT-ICE with LR-iTT-R. Both proposed methods achieve smaller rank than TT-ICE  between the physical-sample interface determining the ROM dimension. Overall, LR-iTT-R offers a favorable balance among reconstruction accuracy, computational efficiency, and compression effectiveness.}

\begin{table}[htbp]
\centering
\caption{\ssy{
Comparison of TT-ICE, LR-iTT-F, and
LR-iTT-R with $\alpha=0.1$ on the homogeneous 2D2V training set under three TT conversion tolerance settings. The mean relative errors are evaluated over the training snapshots. For each setting, the speedup
is the ratio of the average TT-ICE wall time to the average wall time
of the corresponding method.
}}
\label{tab:two-d-two-v-results}
\scriptsize
\setlength{\tabcolsep}{4pt}
\begin{tabular*}{\linewidth}{@{}c@{\extracolsep{\fill}}lcccc@{}}
\toprule
$(\epsilon,\epsilon_{\mathrm{tt}},r_{\max})$
& Method
& TT ranks
& \shortstack{Mean err.}
& Time (s)
& Speedup \\
\midrule
\multirow{3}{*}{$(10^{-4},10^{-5},30)$}
& TT-ICE
& $(11,11,11,6)$
& $1.509 \times 10^{-4}$
& $3.809 \times 10^{1}$
& $1.00 \times 10^{0}$ \\
& LR-iTT-F
& $(14,16,13,4)$
& $3.742 \times 10^{-5}$
& $1.563 \times 10^{-1}$
& $2.44 \times 10^{2}$ \\
& LR-iTT-R
& $(15,16,13,4)$
& $3.903 \times 10^{-5}$
& $1.388 \times 10^{-1}$
& $2.74 \times 10^{2}$ \\
\midrule
\multirow{3}{*}{$(10^{-5},10^{-6},30)$}
& TT-ICE
& $(18,24,22,9)$
& $1.566 \times 10^{-5}$
& $4.185 \times 10^{1}$
& $1.00 \times 10^{0}$ \\
& LR-iTT-F
& $(20,27,24,5)$
& $6.683 \times 10^{-6}$
& $4.165 \times 10^{-1}$
& $1.00 \times 10^{2}$ \\
&LR-iTT-R
& $(20,28,22,5)$
& $7.092 \times 10^{-6}$
& $3.504 \times 10^{-1}$
& $1.19 \times 10^{2}$ \\
\midrule
\multirow{3}{*}{$(10^{-6},10^{-7},60)$}
& TT-ICE
& $(26,47,51,19)$
& $2.943 \times 10^{-6}$
& $4.402 \times 10^{1}$
& $1.00 \times 10^{0}$ \\
& LR-iTT-F
& $(30,50,49,8)$
& $7.510 \times 10^{-7}$
& $2.318 \times 10^{0}$
& $1.90 \times 10^{1}$ \\
& LR-iTT-R
& $(32,54,49,8)$
& $7.110 \times 10^{-7}$
& $2.226 \times 10^{0}$
& $1.98 \times 10^{1}$ \\
\bottomrule
\end{tabular*}
\end{table}

\begin{table}[htbp]
\centering
\caption{\ssy{
Comparison of TT-ICE, LR-iTT-F, and
LR-iTT-R with $\alpha=0.1$ on the homogeneous 2D2V training set without roundings. TT-ICE uses the dense snapshots. Wall times are averaged over three runs.
}}
\label{tab:two-d-two-v-without-truncation}
\scriptsize
\setlength{\tabcolsep}{4pt}
\begin{tabular*}{\linewidth}{@{}c@{\extracolsep{\fill}}lcccc@{}}
\toprule
$(\epsilon,\epsilon_{\mathrm{tt}})$
& Method
& TT ranks
& \shortstack{Mean err.}
& Time (s)
& Speedup \\
\midrule
\multirow{3}{*}{$(10^{-4},10^{-5})$}
& TT-ICE
& $(11,11,11,6)$
& $1.509 \times 10^{-4}$
& $3.944 \times 10^{1}$
& $1.00 \times 10^{0}$ \\
& LR-iTT-F
& $(26,42,33,4)$
& $3.283 \times 10^{-5}$
& $7.362 \times 10^{-1}$
& $5.36 \times 10^{1}$ \\
& LR-iTT-R
& $(15,16,13,4)$
& $3.903 \times 10^{-5}$
& $2.316 \times 10^{-1}$
& $1.70 \times 10^{2}$ \\
\midrule
\multirow{3}{*}{$(10^{-5},10^{-6})$}
& TT-ICE
& $(18,24,22,9)$
& $1.566 \times 10^{-5}$
& $4.032 \times 10^{1}$
& $1.00 \times 10^{0}$ \\
& LR-iTT-F
& $(39,73,69,5)$
& $5.445 \times 10^{-6}$
& $1.088 \times 10^{1}$
& $3.71 \times 10^{0}$ \\
& LR-iTT-R
& $(21,31,25,5)$
& $6.667 \times 10^{-6}$
& $3.809 \times 10^{-1}$
& $1.06 \times 10^{2}$ \\
\bottomrule
\end{tabular*}
\end{table}

\subsubsection{Changing-scattering medium}

We consider the changing-scattering dataset with the configuration given in Table~\ref{tab:rte-benchmark-configs}. As in the homogeneous example, the low-rank solver \cite{guo2026highly} generates the solution snapshots in low-rank matrix-factorized form. We then construct four TT versions of this dataset by converting each snapshot to TT format using four conversion tolerances \(10^{-2}\), \(10^{-3}\), \(10^{-4}\), and \(10^{-5}\), respectively. After conversion, each TT snapshot has mode size $n=(100,100,16,32)$
and TT rank structure $(1,r_1,r_2,r_3,1)$. Compared with the homogeneous example, this dataset exhibits substantially larger TT ranks. The corresponding mean internal TT ranks $(r_1,r_2,r_3)$ of the four TT datasets are
\[
\begin{aligned}
&(6.00,\,25.39,\,14.57), \qquad
(10.00,\,51.24,\,22.24),\\
&(15.00,\,86.35,\,28.51), \qquad
(21.00,\,116.94,\,32.00).
\end{aligned}
\]

For each of the four data sets, all three methods use the same nominal
relative tolerance $\epsilon$, TT-rounding tolerance
$\epsilon_{\mathrm{tt}}$, and rank trigger $r_{\max}$ listed in
Table~\ref{tab:changing-scattering-results}. Intermediate TT-rounding
is applied when triggered, while the final TT-rounding step is omitted.

\textbf{Reconstruction accuracy, TT ranks, and wall time.}
\ssy{Table~\ref{tab:changing-scattering-results} summarizes the mean reconstruction error, TT ranks of the accumulated tensor, and wall time for the four TT datasets. Both proposed methods achieve minorly lower mean errors and shorter wall times than TT-ICE in all four settings. The speedups of LR-iTT-F are $172$, $43.4$, $6.56$, and $2.24$, respectively; those of LR-iTT-R are $182$, $47.4$, $6.83$, and $2.37$. At the loosest setting, the mean errors of LR-iTT-F and LR-iTT-R are $5.567\times10^{-3}$ and $5.591\times10^{-3}$, compared with $1.323\times10^{-2}$ for TT-ICE. At the most stringent setting, LR-iTT-R takes $14.29$ seconds versus $33.85$ seconds for TT-ICE and achieves a lower mean error ($9.701\times10^{-6}$ versus $1.676\times10^{-5}$). LR-iTT-R has a slightly lower observed mean wall time than LR-iTT-F in all four settings, with mean errors of similar magnitude.}

\begin{figure}[t]
\centering
\includegraphics[width=\textwidth]{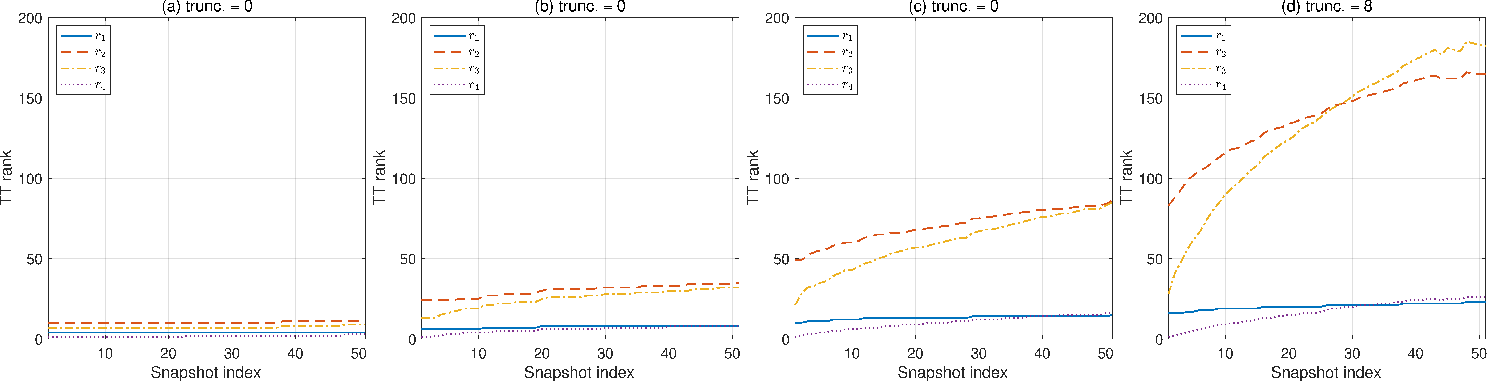}
\par\vspace{4pt}
\includegraphics[width=\textwidth]{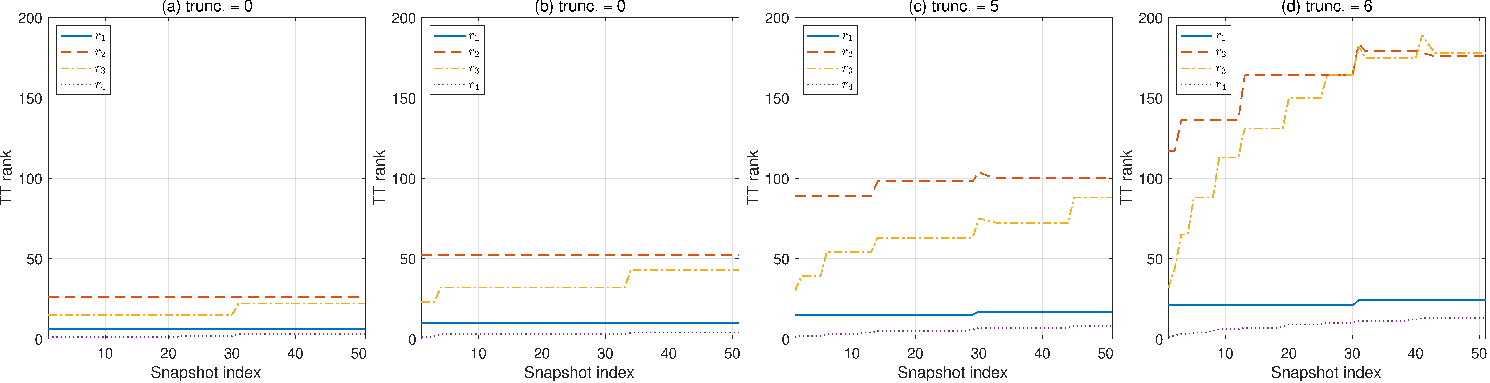}
\par\vspace{4pt}
\includegraphics[width=\textwidth]{Figuras/ssy_plots/fig6_ttice.eps}
\caption{\ssy{Evolution of the four internal TT ranks on the changing-scattering
2D2V training set under the three parameter settings in
Table~\ref{tab:changing-scattering-results}. The top, middle, and bottom rows
correspond to LR-iTT-F, LR-iTT-R, and TT-ICE, respectively.
In each row, the panels correspond to:
(a) $\epsilon=10^{-2}$, $\epsilon_{\rm tt}=10^{-3}$, $r_{\max}=60$;
(b) $\epsilon=10^{-3}$, $\epsilon_{\rm tt}=10^{-4}$, $r_{\max}=60$;
(c) $\epsilon=10^{-4}$, $\epsilon_{\rm tt}=10^{-5}$, $r_{\max}=100$;
(d) $\epsilon=10^{-5}$, $\epsilon_{\rm tt}=10^{-6}$, $r_{\max}=180$.}}
\label{fig:changing-scattering-evo}
\end{figure}

\ssy{The rank behavior is also consistent with that observed in the homogeneous example. TT-ICE generally has smaller ranks at the interfaces between modes. Its ranks at the physical-sample interface are equal to or larger than those of both proposed methods. At the most stringent setting, this interface rank is $26$ for TT-ICE and $13$ for both LR-iTT-F and LR-iTT-R. As the tolerances are tightened, the accumulated ranks and the cost of core-wise projections increase, and the observed speedups decrease. Even at the most stringent setting, LR-iTT-F and LR-iTT-R remain more than 2 times faster than TT-ICE. Figures~\ref{fig:changing-scattering-evo} show how these internal ranks evolve as the snapshots are processed.}


\begin{table}[htbp]
\centering
\caption{\ssy{
Comparison of TT-ICE, LR-iTT-F, and
LR-iTT-R with $\alpha=0.1$ on the changing-scattering 2D2V training set
under four TT conversion tolerance settings. The mean relative errors
are evaluated over the training snapshots. For each setting, the speedup
is the ratio of the average TT-ICE wall time to the average wall time
of the corresponding method.
}}
\label{tab:changing-scattering-results}
\scriptsize
\setlength{\tabcolsep}{4pt}
\begin{tabular*}{\linewidth}{@{}c@{\extracolsep{\fill}}lcccc@{}}
\toprule
$(\epsilon,\epsilon_{\mathrm{tt}},r_{\max})$
& Method
& TT ranks
& \shortstack{Mean rel.\\error}
& Time (s)
& Speedup \\
\midrule
\multirow{3}{*}{$(10^{-2},10^{-3},60)$}
& TT-ICE
& $(4,11,9,3)$
& $1.323 \times 10^{-2}$
& $2.463 \times 10^{1}$
& $1.00 \times 10^{0}$ \\
& LR-iTT-F
& $(6,26,29,3)$
& $5.567 \times 10^{-3}$
& $1.433 \times 10^{-1}$
& $1.72 \times 10^{2}$ \\
& LR-iTT-R
& $(6,26,22,3)$
& $5.591 \times 10^{-3}$
& $1.351 \times 10^{-1}$
& $1.82 \times 10^{2}$ \\
\midrule
\multirow{3}{*}{$(10^{-3},10^{-4},60)$}
& TT-ICE
& $(8,35,32,8)$
& $1.601 \times 10^{-3}$
& $2.604 \times 10^{1}$
& $1.00 \times 10^{0}$ \\
& LR-iTT-F
& $(10,51,41,4)$
& $8.029 \times 10^{-4}$
& $5.996 \times 10^{-1}$
& $4.34 \times 10^{1}$ \\
& LR-iTT-R
& $(10,52,43,4)$
& $7.929 \times 10^{-4}$
& $5.496 \times 10^{-1}$
& $4.74 \times 10^{1}$ \\
\midrule
\multirow{3}{*}{$(10^{-4},10^{-5},100)$}
& TT-ICE
& $(15,86,85,16)$
& $1.685 \times 10^{-4}$
& $2.978 \times 10^{1}$
& $1.00 \times 10^{0}$ \\
& LR-iTT-F
& $(17,100,80,6)$
& $8.652 \times 10^{-5}$
& $4.538 \times 10^{0}$
& $6.56 \times 10^{0}$ \\
& LR-iTT-R
& $(17,100,88,8)$
& $9.178 \times 10^{-5}$
& $4.359 \times 10^{0}$
& $6.83 \times 10^{0}$ \\
\midrule
\multirow{3}{*}{$(10^{-5},10^{-6},180)$}
& TT-ICE
& $(23,165,182,26)$
& $1.676 \times 10^{-5}$
& $3.385 \times 10^{1}$
& $1.00 \times 10^{0}$ \\
& LR-iTT-F
& $(24,174,176,13)$
& $9.487 \times 10^{-6}$
& $1.514 \times 10^{1}$
& $2.24 \times 10^{0}$ \\
& LR-iTT-R
& $(24,176,178,13)$
& $9.701 \times 10^{-6}$
& $1.429 \times 10^{1}$
& $2.37 \times 10^{0}$ \\
\bottomrule
\end{tabular*}
\end{table}

\subsection{Applications in the ROM\label{sec:numerical-rom}}
We next apply the proposed incremental TT compression to build reduced bases for
ROMs, following Section~\ref{sec:tt_rom}, on the three RTE benchmarks in
Section~\ref{sec:benchmark}. For each example, 
the low-rank solver in \cite{guo2026highly} generates
the angular-flux snapshots in the form of low-rank matrix factorization:
\begin{equation}
F_j=U_j\operatorname{diag}(s_j)V_j^T\in\mathbb{R}^{(4N_xN_y)\times (N_\theta N_{v_z})},\; U_j\in\mathbb{R}^{4N_xN_y\times r_j},\;V_j\in\mathbb{R}^{N_\theta N_{v_z}\times r_j},
\end{equation}
where $U_j$ and $V_j$ are low-rank factors collecting spatial and angular DOFs, respectively.
In the
projection-based ROM, these low-rank matrix factors are used directly as
streaming input. In interpolation-based ROM, the same snapshots are
converted to order-four TT tensors without additional truncation to demonstrate the flexibility of the
proposed construction. In both cases, the online tests use five randomly
sampled parameters outside the training set.

\subsubsection{Applications in the projection-based ROM\label{sec:projection-rom}}
We begin with the projection-based intrusive setting, where the reduced basis is used to project
the affine operators, as well as the source term.

\textbf{Incremental ROM construction from low-rank data.}
We construct the reduced basis directly from the streaming low-rank snapshots
with residual enrichment threshold \(\epsilon=10^{-5}\), without prescribing a
maximum TT rank.  This setting highlights two advantages.
\begin{enumerate}
\item \textbf{Low-rank data generation.}
Table~\ref{tab:projection-lowrank-data} shows that the low-rank solver is
faster in all three examples, with speedups ranging from \(1.46\) for the
lattice case to \(6.59\) for the homogeneous case, and reduces the
storage by factors between \(2.00\) and \(7.41\).  The compressed low-rank
solutions remain close to the full-rank references: The mean relative errors are of order $10^{-8}$--$10^{-6}$
for the angular flux $f$ and $10^{-8}$--$10^{-7}$
for the scalar flux $\rho$.

\item \textbf{Incremental TT compression of the training data.}
\ssy{Table~\ref{tab:projection-rom-basis} summarizes the reduced bases constructed by LR-iTT-F and LR-iTT-R with $\alpha=0.1$. For the homogeneous, changing-scattering, and lattice examples, LR-iTT-F yields accumulated TT ranks $(1,122,3,1)$, $(1,372,14,1)$, and $(1,1130,21,1)$,
with ROM dimensions $3$, $14$, and $21$, respectively.
LR-iTT-R yields TT ranks $(1,69,4,1)$, $(1,272,14,1)$, and $(1,837,25,1)$, with ROM dimensions $4$, $14$, and $25$.
Despite these equal or slightly larger ROM dimensions, LR-iTT-R increases the storage reduction factors from
$157.31$, $40.88$, and $21.72$ to
$270.39$, $55.91$, and $27.75$, respectively.}
\end{enumerate}

After constructing the TT reduced basis $U_{\rm rb}$, the affine linear operators and source term are projected using the left-to-right TT contractions described in Section~\ref{sec:projection-rom}.

\begin{table}[t]
\caption{\ssy{Average performance of low-rank data generation over five online test parameters.  The errors are measured relative to the corresponding full-rank solutions.  The storage reduction is the ratio of dense-matrix
storage to the storage of the modal low-rank factors.}}
\label{tab:projection-lowrank-data}
\centering
\small
\resizebox{\textwidth}{!}{%
\begin{tabular}{lccccccc}
\toprule
Example & Full solve time (s) & Low-rank solve time (s) &
Speedup & Mean rank & Storage reduction &
Low-rank \(f\) error & Low-rank \(\rho\) error \\
\midrule
Homogeneous & \(4.140\times 10^{2}\) & \(6.283\times 10^{1}\) &
\(6.59\) & \(67.0\) & \(7.41\) &
\(6.887\times 10^{-8}\) & \(6.744\times 10^{-8}\) \\
Changing scattering & \(2.377\times 10^{2}\) & \(6.739\times 10^{1}\) &
\(3.53\) & \(120.4\) & \(4.04\) &
\(4.136\times 10^{-6}\) & \(2.992\times 10^{-7}\) \\
Lattice & \(4.055\times 10^{2}\) & \(2.769\times 10^{2}\) &
\(1.46\) & \(141.2\) & \(2.00\) &
\(1.384\times 10^{-6}\) & \(3.013\times 10^{-7}\) \\
\bottomrule
\end{tabular}
}
\end{table}


\begin{table}[htbp]
    \centering
    \caption{\ssy{
        Offline construction of the projection-based ROM using
        LR-iTT-F and LR-iTT-R with $\alpha=0.1$. Both methods use residual
        tolerance $\epsilon=10^{-5}$. The TT ranks are those of the matrix-TT reduced basis $U_{\mathrm{rb}}$, formed from the first two cores of the accumulated tensor. The storage reductions compare the full-rank training snapshot collection with the accumulated TT ROM representation $\psi$.
    }}
    \label{tab:projection-rom-basis}

    \small
    \setlength{\tabcolsep}{3pt}
    \resizebox{\linewidth}{!}{%
    \begin{tabular}{@{}lccclccc@{}}
        \toprule
        Example
        & \shortstack{Snapshots}
        & \shortstack{Matrix size}
        & \shortstack{Mean input rank}
        & Method
        & \shortstack{Ranks of $U_{\mathrm{rb}}$}
        & \shortstack{ROM dim}
        & \shortstack{Storage reduction} \\
        \midrule

        \multirow{2}{*}{Homogeneous}
        & \multirow{2}{*}{$41$}
        & \multirow{2}{*}{$16384 \times 512$}
        & \multirow{2}{*}{$65.37$}
        & LR-iTT-F 
        & $(1,122,3)$
        & $3$
        & $157.31$ \\
        & & & &
        LR-iTT-R
        & $(1, 69, 4)$
        & $4$
        & $270.39$ \\

        \midrule

        \multirow{2}{*}{Changing scattering}
        & \multirow{2}{*}{$51$}
        & \multirow{2}{*}{$10000 \times 512$}
        & \multirow{2}{*}{$121.53$}
        & LR-iTT-F
        & $(1,372,14)$
        & $14$
        & $40.88$ \\
        & & & &
        LR-iTT-R
        & $(1, 272, 14)$
        & $14$
        & $55.91$ \\

        \midrule

        \multirow{2}{*}{Lattice}
        & \multirow{2}{*}{$121$}
        & \multirow{2}{*}{$14400 \times 288$}
        & \multirow{2}{*}{$140.76$}
        & LR-iTT-F
        & $(1,1130,21)$
        & $21$
        & $21.72$ \\
        & & & &
        LR-iTT-R
        & $(1, 837, 25)$
        & $25$
        & $27.75$ \\

        \bottomrule
    \end{tabular}%
    }
\end{table}

\textbf{Online efficiency and accuracy.}
\ssy{Table~\ref{tab:projection-rom-online} summarizes the online performance of the projection-based ROM. The online stage consists of solving the reduced system
and reconstructing the solution from the reduced coefficients. After model reduction, the reconstruction step becomes the dominant online cost.
Nevertheless, the overall online speedups for LR-iTT-F are substantial:
\(1.85\times10^5\), \(1.99\times10^4\), and \(1.16\times10^4\) over the
full-rank solver, and \(2.81\times10^4\), \(5.61\times10^3\), and
\(7.84\times10^3\) over the low-rank solver for the homogeneous,
changing-scattering, and lattice examples, respectively. LR-iTT-R further reduces this reconstruction time in each case. Its speedups over the full-rank solver are
$2.01\times10^{5}$, $2.84\times10^{4}$, and $1.65\times10^{4}$
for three examples, respectively. The corresponding speedups over the low-rank solver are
$3.06\times10^{4}$, $8.02\times10^{3}$, and $1.13\times10^{4}$. This substantial acceleration is achieved without loss of accuracy: the mean
relative errors of both approaches retain order $10^{-6}$--$10^{-5}$ for both $f$ and $\rho$.}


\begin{table}[htbp]
    \centering
    \caption{\ssy{
        Average online performance of the projection-based ROM
        over five test parameters. The reduced bases are constructed using LR-iTT-F and LR-iTT-R with $\alpha=0.1$.
        Total online time includes the reduced solve
        and scalar-flux reconstruction. Speedups are the average ratios of the full-rank or
        low-rank solve time to the total online time.
        The mean relative errors are measured against the corresponding full-rank solutions.
    }}
    \label{tab:projection-rom-online}

    \small
    \setlength{\tabcolsep}{4pt}
    \resizebox{\linewidth}{!}{%
    \begin{tabular}{@{}llcccccc@{}}
        \toprule
        Example
        & Method
        & \shortstack{Reduced solve\\time (s)}
        & \shortstack{Scalar-flux\\reconstruction (s)}
        & \shortstack{Speedup over\\full-rank solver}
        & \shortstack{Speedup over\\low-rank solver}
        & \shortstack{Mean rel.\\error in $f$}
        & \shortstack{Mean rel.\\error in $\rho$} \\
        \midrule

        \multirow{2}{*}{Homogeneous}
        & LR-iTT-F
        & $4.733 \times 10^{-4}$
        & $2.032 \times 10^{-3}$
        & $1.852 \times 10^{5}$
        & $2.809 \times 10^{4}$
        & $1.612 \times 10^{-5}$
        & $1.612 \times 10^{-5}$ \\
        & LR-iTT-R
        & $1.061 \times 10^{-3}$
        & $1.694 \times 10^{-3}$
        & $2.015 \times 10^{5}$
        & $3.058 \times 10^{4}$
        & $1.462 \times 10^{-5}$
        & $1.460 \times 10^{-5}$ \\

        \midrule

        \multirow{2}{*}{Changing scattering}
        & LR-iTT-F
        & $6.856 \times 10^{-4}$
        & $1.236 \times 10^{-2}$
        & $1.988 \times 10^{4}$
        & $5.613 \times 10^{3}$
        & $4.886 \times 10^{-6}$
        & $2.910 \times 10^{-6}$ \\
        & LR-iTT-R
        & $4.080 \times 10^{-4}$
        & $8.373 \times 10^{-3}$
        & $2.844 \times 10^{4}$
        & $8.020 \times 10^{3}$
        & $5.151 \times 10^{-6}$
        & $3.052 \times 10^{-6}$ \\

        \midrule

        \multirow{2}{*}{Lattice}
        & LR-iTT-F
        & $1.258 \times 10^{-4}$
        & $3.561 \times 10^{-2}$
        & $1.157 \times 10^{4}$
        & $7.844 \times 10^{3}$
        & $2.702 \times 10^{-6}$
        & $1.495 \times 10^{-6}$ \\
        & LR-iTT-R
        & $1.199 \times 10^{-4}$
        & $2.511 \times 10^{-2}$
        & $1.648 \times 10^{4}$
        & $1.133 \times 10^{4}$
        & $4.072 \times 10^{-6}$
        & $2.061 \times 10^{-6}$ \\

        \bottomrule
    \end{tabular}%
    }
\end{table}

\subsubsection{Applications in the interpolation-based ROM\label{sec:interpolation-rom}}
To demonstrate the flexibility of our method, we consider the non-intrusive interpolation-based ROM using TT snapshots with nodal modes \((2N_x,2N_y,N_{v_z},N_\theta)\).
The accumulated sample-augmented TT tensor provides
both the reduced basis and the training coefficient core, which is interpolated
online for new parameters.

\textbf{Offline TT basis construction.}
\ssy{Table~\ref{tab:interpolation-rom-basis} summarizes the offline construction. For the homogeneous, changing-scattering, and lattice examples, respectively,
LR-iTT-F produces ROM dimensions of $3$, $13$, and $20$, with storage reduction factors of $274.80$, $69.54$, and $22.79$. LR-iTT-R with $\alpha=0.1$ produces ROM dimensions of $5$, $14$, and $24$, while increasing the corresponding storage reduction factors to $526.20$,
$127.19$, and $29.39$. Thus, despite the larger ROM dimensions, LR-iTT-R
achieves a more compact representation of the
accumulated TT tensor.}

\begin{table}[htbp]
    \centering
    \caption{\ssy{
        Offline construction of the interpolation-based ROM using
        LR-iTT-F and LR-iTT-R with $\alpha=0.1$. Both methods use residual tolerance $\epsilon=10^{-5}$
without TT rounding. The mean TT input ranks are computed after converting the nodal low-rank matrix snapshots into order-four TT tensors. The storage reductions compare the full-rank training
        snapshot collection with the accumulated TT ROM representation $\psi$.
    }}
    \label{tab:interpolation-rom-basis}

    \small
    \setlength{\tabcolsep}{3pt}
    \resizebox{\linewidth}{!}{%
    \begin{tabular}{@{}lccclccc@{}}
        \toprule
        Example
        & \shortstack{Snapshots}
        & \shortstack{TT modes}
        & \shortstack{Mean input ranks}
        & Method
        & \shortstack{Ranks of $U_{\mathrm{rb}}$}
        & \shortstack{ROM dim}
        & \shortstack{Storage reduction} \\
        \midrule

        \multirow{2}{*}{Homogeneous}
        & \multirow{2}{*}{$41$}
        & \multirow{2}{*}{$128 \times 128 \times 16 \times 32$}
        & \multirow{2}{*}{$(1,67.71,65.37,32,1)$}
        & LR-iTT-F
        & $(1,67,122,96,3)$
        & $3$
        & $274.80$ \\
        & & & &
        LR-iTT-R
        & $(1, 67, 69, 42, 5)$
        & $5$
        & $526.20$\\

        \midrule

        \multirow{2}{*}{Changing scattering}
        & \multirow{2}{*}{$51$}
        & \multirow{2}{*}{$100 \times 100 \times 16 \times 32$}
        & \multirow{2}{*}{$(1,52.41,121.53,32,1)$}
        & LR-iTT-F
        & $(1,55,359,288,13)$
        & $13$
        & $69.54$ \\
        & & & &
        LR-iTT-R
        & $(1, 55, 213, 227, 14 )$
        & $14$
        & $127.19$ \\

        \midrule

        \multirow{2}{*}{Lattice}
        & \multirow{2}{*}{$121$}
        & \multirow{2}{*}{$120 \times 120 \times 12 \times 24$}
        & \multirow{2}{*}{$(1,120,140.76,24,1)$}
        & LR-iTT-F
        & $(1,120,1130,408,20)$
        & $20$
        & $22.79$ \\
        & & & &
        LR-iTT-R
        & $(1, 120, 835, 475, 24)$
        & $24$
        & $29.39$\\

        \bottomrule
    \end{tabular}%
    }
\end{table}

\textbf{Online coefficient interpolation and output reconstruction.}
We use the same randomly sampled test parameters as in
Section~\ref{sec:projection-rom}. Following Section~\ref{sec:non-intrusive}, the coefficient core is
interpolated using cubic splines for the one-parameter homogeneous and
changing-scattering examples, and tensor-product cubic splines in
\((\sigma_s,\sigma_a)\) for the lattice example. The angular flux and scalar
flux are then recovered from the compressed TT representation by contractions
with the TT basis and the angular quadrature weights, respectively.

\ssy{Table~\ref{tab:interpolation-rom-online} reports the average online cost and accuracy over five test parameters. For both methods, computing the interpolation weights takes substantially less time than evaluating the scalar flux from the reduced TT representation. For the homogeneous, changing-scattering, and lattice examples,
respectively, LR-iTT-F achieves speedups of
$2.67\times10^{4}$, $1.26\times10^{4}$, and $4.32\times10^{3}$
over the full-rank solver, and
$4.06\times10^{3}$, $3.58\times10^{3}$, and $2.94\times10^{3}$
over the low-rank solver. With $\alpha=0.1$, LR-iTT-R reduces the measured scalar-flux
evaluation time in all three examples, increasing the corresponding
speedups to
$2.96\times10^{4}$, $1.88\times10^{4}$, and $5.18\times10^{3}$
over the full-rank solver, and
$4.51\times10^{3}$, $5.32\times10^{3}$, and $3.53\times10^{3}$
over the low-rank solver. Both methods retain small mean relative errors
of order $10^{-6}$--$10^{-5}$ for both $f$ and $\rho$.}

\ssy{\textbf{Sensitivity to $\alpha$.}
Tables~\ref{tab:alpha-homogeneous} and \ref{tab:alpha-changing}
assess the influence of $\alpha$ on ROM test accuracy, offline
update time, and TT storage using the same five test parameters
as above. Increasing $\alpha$ reduces the offline update time,
but the accuracy and storage do not vary monotonically.
For the changing-scattering example, reducing $\alpha$ from
$0.1$ to $0.05$ lowers the mean errors in $f$ and $\rho$ by
approximately $24\%$ and $38\%$, respectively, while increasing
both update time and storage by approximately $11\%$.
Further reducing $\alpha$ to $0.025$ increases the costs without
improving accuracy. For $\alpha>0.1$, update time decreases
further and the error in $f$ increases, whereas the error in
$\rho$ remains non-monotone. Storage generally decreases,
with a local increase at $\alpha=0.2$.}

\pzc{
At $\alpha=0.1$, LR-iTT-R retains accuracy comparable to LR-iTT-F
while reducing offline update time by approximately $56\%$
and $32\%$, and storage by approximately $48\%$ and $45\%$,
for the homogeneous and changing-scattering examples,
respectively. These results support $\alpha=0.1$ as a practical
common choice that balances test accuracy, offline computational
cost, and storage across the two examples.}


\begin{table}[htbp]
    \centering
    \caption{\ssy{
        Average online performance of the interpolation-based ROM over five test parameters. The reduced bases are constructed using LR-iTT-F and LR-iTT-R with $\alpha=0.1$. The interpolation time measures the computation of parameter interpolation weights. The scalar-flux evaluation time includes the TT contractions required to compute $\rho$ and evaluations for $\rho$.
        Speedups are the average ratios of the full-rank or low-rank solve time to the sum of these
        two online costs. The mean relative errors are measured against the reference nodal test snapshots.
    }}
    \label{tab:interpolation-rom-online}

    \small
    \setlength{\tabcolsep}{4pt}
    \resizebox{\linewidth}{!}{%
    \begin{tabular}{@{}llcccccc@{}}
        \toprule
        Example
        & Method
        & \shortstack{Interpolation\\time (s)}
        & \shortstack{Scalar-flux\\evaluation (s)}
        & \shortstack{Speedup over\\full-rank solver}
        & \shortstack{Speedup over\\low-rank solver}
        & \shortstack{Mean rel.\\error in $f$}
        & \shortstack{Mean rel.\\error in $\rho$} \\
        \midrule

        \multirow{2}{*}{Homogeneous}
        & LR-iTT-F
        & $5.185 \times 10^{-4}$
        & $1.593 \times 10^{-2}$
        & $2.674 \times 10^{4}$
        & $4.061 \times 10^{3}$
        & $5.109 \times 10^{-6}$
        & $3.985 \times 10^{-6}$ \\
        & LR-iTT-R
        & $4.802 \times 10^{-4}$
        & $1.385 \times 10^{-2}$
        & $2.963 \times 10^{4}$
        & $4.506 \times 10^{3}$
        & $7.981 \times 10^{-6}$
        & $2.823 \times  10^{-6}$ \\

        \midrule

        \multirow{2}{*}{Changing scattering}
        & LR-iTT-F
        & $5.719 \times 10^{-4}$
        & $1.832 \times 10^{-2}$
        & $1.264 \times 10^{4}$
        & $3.583 \times 10^{3}$
        & $8.072 \times 10^{-6}$
        & $4.652 \times 10^{-6}$ \\
        & LR-iTT-R
        & $6.301 \times 10^{-4}$
        & $1.204 \times 10^{-2}$
        & $1.878 \times 10^{4}$
        & $5.317 \times 10^{3}$
        & $8.311 \times 10^{-6}$
        & $3.034 \times 10^{-6}$ \\

        \midrule

        \multirow{2}{*}{Lattice}
        & LR-iTT-F
        & $1.358 \times 10^{-3}$
        & $9.370 \times 10^{-2}$
        & $4.317 \times 10^{3}$
        & $2.945 \times 10^{3}$
        & $1.080 \times 10^{-5}$
        & $2.982 \times 10^{-6}$ \\
        & LR-iTT-R
        & $1.337 \times 10^{-3}$
        & $7.813 \times 10^{-2}$
        & $5.178 \times 10^{3}$
        & $3.531 \times 10^{3}$
        & $1.084 \times 10^{-5}$
        & $2.346 \times 10^{-6}$ \\

        \bottomrule
    \end{tabular}%
    }
\end{table}

\begin{table}[htbp]
\centering
\caption{\ssy{
Effect of the enrichment parameter $\alpha$ on the interpolation-based ROM for the homogeneous example. The mean relative errors in angular flux $f$ and scalar flux $\rho$ are averaged over five test parameters. Offline update time is the median of three runs over the full training snapshot stream. Storage is the total number of TT-core elements.
}}
\label{tab:alpha-homogeneous}
\scriptsize
\setlength{\tabcolsep}{4pt}
\begin{tabular*}{\linewidth}{@{}l@{\extracolsep{\fill}}ccccc@{}}
    \toprule
    Approach
    & \shortstack{Mean err. $f$}
    & \shortstack{Mean err. $\rho$}
    & \shortstack{Time (s)}
    & TT ranks
    & \shortstack{Storage} \\
    \midrule
    LR-iTT-F
    & $5.109 \times 10^{-6}$ & $3.985 \times 10^{-6}$ & 8.29
    & $(67,122,96,3)$ & 1,251,579 \\
    LR-iTT-R ($\alpha=0.025$)
    & $7.871 \times 10^{-6}$ & $2.791 \times 10^{-6}$ & 3.70
    & $(67,71,47,5)$ & 678,589 \\
    LR-iTT-R ($\alpha=0.05$)
    & $7.912 \times 10^{-6}$ & $2.800 \times 10^{-6}$ & 3.64
    & $(67,69,45,5)$ & 657,405 \\
    LR-iTT-R ($\alpha=0.1$)
    & $7.981 \times 10^{-6}$ & $2.823 \times 10^{-6}$ & 3.62
    & $(67,69,42,5)$ & 653,613 \\
    LR-iTT-R ($\alpha=0.2$)
    & $8.730 \times 10^{-6}$ & $2.823 \times 10^{-6}$ & 3.58
    & $(67,67,39,5)$ & 631,421 \\
    LR-iTT-R ($\alpha=0.5$)
    & $9.481 \times 10^{-6}$ & $4.153 \times 10^{-6}$ & 3.55
    & $(67,66,38,5)$ & 621,005 \\
    LR-iTT-R ($\alpha=1$)
    & $1.242 \times 10^{-5}$ & $5.061 \times 10^{-6}$ & 3.52
    & $(67,66,37,4)$ & 618,564 \\
    \bottomrule
\end{tabular*}
\end{table}

\begin{table}[htbp]
\centering
\caption{\ssy{
Effect of the enrichment parameter $\alpha$ on the interpolation-based ROM for the changing-scattering example. The mean relative errors in angular flux $f$ and scalar flux $\rho$ are averaged over five test parameters. Offline update time is the median of three runs over the full training snapshot stream. Storage is the total number of TT-core elements.
}}
\label{tab:alpha-changing}
\scriptsize
\setlength{\tabcolsep}{4pt}
\begin{tabular*}{\linewidth}{@{}l@{\extracolsep{\fill}}ccccc@{}}
    \toprule
    Approach
    & \shortstack{Mean err. $f$}
    & \shortstack{Mean err. $\rho$}
    & \shortstack{Time (s)}
    & TT ranks
    & \shortstack{Storage} \\
    \midrule
    LR-iTT-F
    & $8.072 \times 10^{-6}$ & $4.652 \times 10^{-6}$ & 20.33
    & $(55,359,288,13)$ & 3,754,743 \\
    LR-iTT-R ($\alpha=0.025$)
    & $8.473 \times 10^{-6}$ & $4.637 \times 10^{-6}$ & 18.48
    & $(55,331,266,13)$ & 3,346,055 \\
    LR-iTT-R ($\alpha=0.05$)
    & $6.330 \times 10^{-6}$ & $1.876 \times 10^{-6}$ & 15.47
    & $(55,231,244,14)$ & 2,287,850 \\
    LR-iTT-R ($\alpha=0.1$)
    & $8.311 \times 10^{-6}$ & $3.034 \times 10^{-6}$ & 13.89
    & $(55,213,227,14)$ & 2,053,026 \\
    LR-iTT-R ($\alpha=0.2$)
    & $9.354 \times 10^{-6}$ & $2.541 \times 10^{-6}$ & 12.41
    & $(55,226,229,14)$ & 2,179,870 \\
    LR-iTT-R ($\alpha=0.5$)
    & $1.137 \times 10^{-5}$ & $2.893 \times 10^{-6}$ & 10.92
    & $(55,186,184,17)$ & 1,677,047 \\
    LR-iTT-R ($\alpha=1$)
    & $1.402 \times 10^{-5}$ & $2.438 \times 10^{-6}$ & 10.15
    & $(55,177,176,23)$ & 1,608,141 \\
    \bottomrule
\end{tabular*}
\end{table}

\section{Conclusions\label{sec:conclusions}}
We developed deterministic incremental TT compression algorithms (LR-iTT-F and LR-iTT-R) for streaming
high-dimensional data that are already available in TT or low-rank matrix factorization
formats. For streaming TT data, the method updates an accumulated TT representation directly at the
core level through projection, residual orthogonalization, and adaptive
enrichment. In this way, it avoids reconstructing both the incoming tensors and
the accumulated full tensor, preserving the computational advantage of low-rank
representations. We also established approximation error bounds and computational complexity estimates for the
proposed approach, and demonstrated how the accumulated
TT representation can be exploited to construct ROMs directly from low-rank data through operations on TT cores.

Numerical experiments on parametric radiative transfer equations confirmed the
accuracy and efficiency of the proposed method. The proposed incremental compression algorithm
achieved reconstruction accuracy comparable to TT-ICE while substantially reducing wall time. In addition, the low-rank ROMs built from the accumulated TT produced accurate reduced solutions with significant
online speedups and storage reductions.

Several directions are worth further exploration. These include extending the
projection-based ROM construction to non-affine or nonlinear operators through
tensor hyper-reduction. Another promising direction is to develop more efficient
dataset-level parallel TT compression algorithms by combining the proposed
incremental TT compression with ideas from TT-ICE.

\subsection*{Acknowledgments}
Part of this research was performed while W. Guo was visiting the Institute for Pure and Applied Mathematics (IPAM), which is supported by the National Science Foundation (Grant No. DMS-1925919/2422832). Z. Peng was partially supported by the Hong Kong Research Grants Council grants Early Career Scheme 26302724 and General Research Fund 16306825. 

\section*{Declaration of generative AI and AI-assisted technologies in the writing process}
During the preparation of this work, the authors used AI-assisted tools to assist with code development, debugging, grammar checking, and readability improvement. All numerical experiments, mathematical formulations, algorithms, proofs, results, and conclusions were reviewed and verified by the authors. The authors assume responsibility for all content.
\appendix
\section{Proof of Lemma \ref{lem:prop}}\label{ap:proof} 
Consider two $d$-dimensional tensors $\bm X,\,\bm Y\in\mathbb{R}^{n_1\times n_2\times\cdots\times n_d}$ represented in TT format given by
\[
\bm X = \mU_1\cdots\mU_{i-1}\mW_{i}^X\cdots\mW_{d}^X,\quad \bm Y = \mU_1\cdots\mU_{i-1}\mW_{i}^Y\cdots\mW_{d}^Y
\]
with TT-ranks $\{r_0,r_1,\ldots,r_{i-1},r_i^X,\ldots, r_{d}^X\}$ and $\{r_0,r_1,\ldots,r_{i-1},r_i^Y,\ldots, r_{d}^Y\}$, respectively.
Assume that the shared first $i-1$ cores $\mU_1,\dots,\mU_{i-1}$ are left-orthogonal. Let $W_i^X$ and $W_i^Y$ denote the left unfoldings of the $i$-th cores $\mW_{i}^X$ and $\mW_{i}^Y$, respectively. If $W_i^X$ and $W_i^Y$ are orthogonal in the sense that
\[
(W_i^X)^T W_i^Y = \mathbf{0},
\]
then we show that the Frobenius inner product of the tensors is zero:
\[
\langle \bm X, \bm Y \rangle_F = 0.
\]

Let the interface matrix
$U_{\le i-1}\in \mathbb{R}^{(n_1\cdots n_{i-1})\times r_{i-1}}$ 
denote the contraction of the shared first $i-1$ cores, and let the interface matrices
\[
W_{\ge i+1}^X \in \mathbb{R}^{(n_{i+1}\cdots n_d)\times r_i^X},
\qquad
W_{\ge i+1}^Y \in \mathbb{R}^{(n_{i+1}\cdots n_d)\times r_i^Y}
\]
denote the contractions of the remaining cores $i+1,\dots,d$ for $\bm X$ and $\bm Y$, respectively. 

Then, the mode-$i$ unfolding of $\bm X$ and $\bm Y$ can be written as
\[
X_{\langle i\rangle}
=
(U_{\le i-1}\otimes I_{n_i})\,W_i^X\,(W_{\ge i+1}^X)^T,
\]
\[
Y_{\langle i\rangle}
=
(U_{\le i-1}\otimes I_{n_i})\,W_i^Y\,(W_{\ge i+1}^Y)^T,
\]
where \(\otimes\) denotes the Kronecker product and $I_s$ denotes the $s$-by-$s$ identity matrix; see \cite{steinlechner2016riemannian}.
Therefore,
\begin{align*}
\langle \bm X,\bm Y\rangle_F
&=
\left\langle X_{\langle i\rangle},Y_{\langle i\rangle}\right\rangle_F \\
&=
\operatorname{tr}\!\left(
\bigl((U_{\le i-1}\otimes I_{n_i})W_i^X (W_{\ge i+1}^X)^T\bigr)^T
\bigl((U_{\le i-1}\otimes I_{n_i})W_i^Y (W_{\ge i+1}^Y)^T\bigr)
\right) \\
&=
\operatorname{tr}\!\left(
W_{\ge i+1}^X
(W_i^X)^T
(U_{\le i-1}\otimes I_{n_i})^T
(U_{\le i-1}\otimes I_{n_i})
W_i^Y
(W_{\ge i+1}^Y)^T
\right).
\end{align*}
Since the first $i-1$ shared cores are left-orthogonal,
\(
U_{\le i-1}^T U_{\le i-1}=I_{r_{i-1}},
\)
\[
(U_{\le i-1}\otimes I_{n_i})^T (U_{\le i-1}\otimes I_{n_i})
=
(U_{\le i-1}^T U_{\le i-1})\otimes I_{n_i}
=
I_{r_{i-1}n_i}.
\]
Thus,
\[
\langle \bm X,\bm Y\rangle_F
=
\operatorname{tr}\!\left(
W_{\ge i+1}^X
(W_i^X)^T
W_i^Y
(W_{\ge i+1}^Y)^T
\right).
\]
Finally, given the orthogonality condition $(W_i^X)^T W_i^Y = 0$, it immediately follows that
$\langle \bm X,\bm Y\rangle_F=0.$

\bibliography{ref_guo.bib,ref_peng.bib,ref_shang.bib}
\bibliographystyle{elsarticle-num} 

\end{document}